\documentclass{article}%
\usepackage{graphicx}
\usepackage{amsmath}
\usepackage{amsfonts}
\usepackage{amssymb}%
\newtheorem{theorem}{Theorem}
{}

\newtheorem{definition}{Definition}

\newtheorem{lemma}{Lemma}
{}

\newtheorem{proposition}{Proposition}
\newtheorem{remark}{Remark}

\newenvironment{proof}[1][Proof]{\textbf{#1.} }{\ \rule{0.5em}{0.5em}}

\begin{document}

\title{Uniform Convergence of the Spectral Expansion for a Differential Operator with
Periodic Matrix Coefficients}
\author{O.A.Veliev\\{\small Depart. of Math, Fac.of Arts and Sci.,}\\{\small Dogus University, Ac\i badem, 34722,}\\{\small Kadik\"{o}y, \ Istanbul, Turkey.}\\\ {\small e-mail: oveliev@dogus.edu.tr}}
\date{}
\maketitle

\begin{abstract}
In this paper, we obtain asymptotic formulas for eigenvalues and
eigenfunctions of the operator generated by a system of ordinary differential
equations with summable coefficients and the quasiperiodic boundary
conditions. Using these asymptotic formulas, we find conditions on the
coefficients for which the root functions of this operator form a Riesz basis.
Then we obtain the uniformly convergent spectral expansion of the differential
operators with the periodic matrix coefficients

\end{abstract}

\section{Introduction}

Let $L(P_{2},P_{3}...,P_{n})$ be the differential operator generated in the
space $L_{2}^{m}(-\infty,\infty)$ by the differential expression
\[
l(y)=y^{(n)}(x)+P_{2}\left(  x\right)  y^{(n-2)}(x)+P_{3}\left(  x\right)
y^{(n-3)}(x)+...+P_{n}(x)y
\]
and $L_{t}(P_{2},P_{3}...,P_{n})$ be the differential operator generated in
$L_{2}^{m}(0,1)$ by the same differential expression and the boundary
conditions
\begin{equation}
U_{\mathbb{\nu}\text{,}t}(y)\equiv y^{(\mathbb{\nu})}\left(  1\right)
-e^{it}y^{(\mathbb{\nu})}\left(  0\right)  =0,\text{ }\mathbb{\nu
}=0,1,...,(n-1),
\end{equation}
where $n\geq2,$ $P_{\mathbb{\nu}}=\left(  p_{\mathbb{\nu},i,j}\right)  $ is a
$m\times m$ matrix with the complex-valued summable entries $p_{\mathbb{\nu
},i,j}$, $P_{\mathbb{\nu}}\left(  x+1\right)  =P_{\mathbb{\nu}}\left(
x\right)  $ for $\mathbb{\nu}=2,3,...n,$ the eigenvalues $\mu_{1},\mu
_{2},...,\mu_{m}$ of the matrix
\[
C=\int_{0}^{1}P_{2}\left(  x\right)  dx
\]
are simple. Here $L_{2}^{m}(a,b)$ is the space of the vector functions
$f=\left(  f_{1},f_{2},...,f_{m}\right)  ,$ where $f_{k}\in L_{2}(a,b)$ for
$k=1,2,...,m,$ with the norm $\left\Vert .\right\Vert $ and inner product
$(.,.)$ defined by%
\[
\left\Vert f\right\Vert ^{2}=\left(  \int_{a}^{b}\left\vert f\left(  x\right)
\right\vert ^{2}dx\right)  ,\text{ }(f,g)=\int_{a}^{b}\left\langle f\left(
x\right)  ,g\left(  x\right)  \right\rangle dx,
\]
where $\left\vert .\right\vert $ and $\left\langle .,.\right\rangle $ are the
norm and inner product in $\mathbb{C}^{m}.$ \ For notational convenience we
identify $L=L(P_{2},P_{3}...,P_{n})$, $L_{t}=L_{t}(P_{2},P_{3}...,P_{n})$ in
the following.

It is well-known that ( see [2, 10] ) the spectrum $\sigma(L)$ of $L$ is the
union of the spectra $\sigma(L_{t})$ of $L_{t}$ for $t\in\lbrack0,2\pi).$ To
construct the uniformly convergent spectral expansion for $L$ we first obtain
the uniform, with respect to $t\in Q_{\varepsilon}(n),$ asymptotic formula for
the eigenvalues and eigenfunctions of $L_{t},$ where
\begin{equation}
Q_{\varepsilon}(2\mu)=\{t\in Q:|t-\pi k|>\varepsilon,\forall k\in
\mathbb{Z\}},\text{ }Q_{\varepsilon}(2\mu+1)=Q,
\end{equation}
$Q$ is compact connected subset of $\mathbb{C}$ containing a neighborhood of
the interval $[-a,2\pi-a]$, $a\in(0,\frac{\pi}{2}),$ $\varepsilon\in
(0,\frac{a}{2})$ and $\mu=1,2,....$ Then we prove that the root functions of
$L_{t}$ for $t\in\mathbb{C}(n)$ form a Riesz basis in $L_{2}^{m}(0,1),$ where
$\mathbb{C}(2\mu)=\mathbb{C}\backslash\{\pi k:k\in\mathbb{Z}\mathbb{\}}$,
$\mathbb{C}(2\mu+1)=\mathbb{C}$.

Let us introduce some preliminary results and describe the scheme of the
paper. Clearly%
\[
\varphi_{k,j,t}(x)=\frac{e_{j}}{\left\Vert e^{itx}\right\Vert }e^{i\left(
2\pi k+t\right)  x}\text{ \ for \ }j=1,2,...,m,
\]
where $e_{1}=(1,0,0,...,0),$ $e_{2}=(0,1,0,...,0),...,e_{m}=(0,0,...,0,1),$
are the normalized eigenfunctions of the operator $L_{t}(0)$ corresponding to
the eigenvalue $\left(  2\pi ki+ti\right)  ^{n}$, where $k\in\mathbb{Z},$ and
the operator $L_{t}(P_{2},...,P_{n})$ is denoted by $L_{t}(0)$ when
$P_{2}(x)=0,...,P_{n}(x)=0.$ It easily follows from the classical
investigations [12, chapter 3, theorem 2] that the boundary conditions (1) are
regular and all large eigenvalues of $L_{t}$ belongs to one of the sequences
\begin{equation}
\{\lambda_{k,1}(t):\mid k\mid\geq N\mathbb{\}},\text{ }\{\lambda_{k,2}(t):\mid
k\mid\geq N\mathbb{\}},...,\text{ }\{\lambda_{k,m}(t):\mid k\mid\geq
N\mathbb{\}}\text{,}%
\end{equation}
where $N\gg1,$ satisfying the following, uniform with respect to $t\in Q,$
asymptotic formulas
\begin{equation}
\lambda_{k,j}(t)=\left(  2\pi ki+ti\right)  ^{n}+O\left(  k^{n-1-\frac{1}{2m}%
}\right)
\end{equation}
for $j=1,2,...,m.$ We say that the formula \ $f(k,t)=O(h(k))$ is uniform with
respect to $t\in Q$ if there exists a positive constant $c$, independent of
$t,$ such that $\mid f(k,t))\mid<c\mid h(k)\mid$ for all $t\in Q$ and $\mid
k\mid\gg1.$

The method proposed here allows us to obtain the asymptotic formulas of high
accuracy for the eigenvalues $\lambda_{k,j}(t)$ and the corresponding
normalized eigenfunctions $\Psi_{k,j,t}(x)$ of $L_{t}$ when $p_{\nu,i,j}\in
L_{1}[0,1]$ for all $\nu,i,j$ . Note that to obtain the asymptotic formulas of
high accuracy by the classical methods it is required \ that $P_{2},$
$P_{3},...,P_{n}$ be differentiable (see [12]). To obtain the asymptotic
formulas for $L_{t}$ we take the operator $L_{t}(C),$ where $L_{t}%
(P_{2},...,P_{n})$ is denoted by $L_{t}(C)$ when $P_{2}(x)=C,$ $P_{3}%
(x)=0,...,P_{n}(x)=0,$ for an unperturbed operator and $L_{t}-L_{t}(C)$ for a
perturbation. One can easily verify that the eigenvalues and normalized
eigenfunctions of $L_{t}(C)$ are%
\begin{equation}
\mu_{k,j}(t)=\left(  2\pi ki+ti\right)  ^{n}+\mu_{j}\left(  2\pi ki+ti\right)
^{n-2}\text{, \ }\Phi_{k,j,t}(x)=\frac{v_{j}}{\left\Vert e^{itx}\right\Vert
}e^{i\left(  2\pi k+t\right)  x}%
\end{equation}
for $k\in\mathbb{Z},$ $j=1,2,...,m,$ where $v_{1},v_{2},...,v_{m}$ are the
normalized eigenvectors of the matrix $C$ corresponding to the eigenvalues
$\mu_{1},\mu_{2},...,\mu_{m}$ respectively.

In section 2 we investigate the operator $L_{t}$ and prove the following 2 theorems.

\begin{theorem}
There exist positive constants $c_{1}$, $N_{0},$ independent of $t,$ such that
if $t\in Q_{\varepsilon}(n)$ and $\mid k\mid\geq N_{0}$ then the following
assertions hold:

$(a)$ The eigenvalue $\lambda_{k,j}(t)$ of $L_{t},$ satisfying (4), lie in
\[%
{\textstyle\bigcup\limits_{s=1,2,...,m}}
(U(\mu_{k,s}(t),c_{1}|k|^{n-3}\ln|k|)),
\]
where $U(\mu,c)=\{\lambda\in\mathbb{C}$: $\mid\lambda-\mu\mid<c\}.$

$(b)$ If $\lambda_{k,j}(t)\in U(\mu_{k,p(j)}(t),c_{1}|k|^{n-3}\ln|k|),$ then
there exists unique eigenfunction $\Psi_{k,j,t}(x)$ corresponding to
$\lambda_{k,j}(t)$ and this eigenfunction satisfies
\begin{equation}
\sup_{x\in\lbrack0,1]}\mid\Psi_{k,j,t}(x)-\frac{v_{p(j)}}{\left\Vert
e^{itx}\right\Vert }e^{i\left(  2\pi k+t\right)  x}\mid\leq\frac{c_{2}\ln
|k|}{|k|},
\end{equation}
where $c_{2}$ is a constant independent of $t$ and $j.$
\end{theorem}

Note that here and in forthcoming relations we denote by $c_{i}$ for
$i=1,2,...$, the positive constants, independent of $t,$ whose exact values
are inessential. Using Theorem 1 and investigating associated functions of
$L_{t}$ we prove:

\begin{theorem}
$(a)$ The large eigenvalues of $L_{t}$ consist of $m$ sequences (3) satisfying
the following, uniform with respect to $t\in Q_{\varepsilon}(n),$ formula
\begin{equation}
\lambda_{k,j}(t)=\left(  2\pi ki+ti\right)  ^{n}+\mu_{j}\left(  2\pi
ki+ti\right)  ^{n-2}+O(k^{n-3}\ln|k|),
\end{equation}
namely, $\lambda_{k,j}(t)\in U(\mu_{k,j}(t),c_{1}|k|^{n-3}\ln|k|)$ for $\mid
k\mid\geq N_{0},$ where $c_{1}$ and $N_{0}$ are defined in Theorem 1. If $\mid
k\mid\geq N_{0},$ then $\lambda_{k,j}(t)$ for $t\in Q_{\varepsilon}(n)$ is a
simple eigenvalue of $L_{t}$ and the corresponding normalized eigenfunction
$\Psi_{k,j,t}(x)$ satisfies
\begin{equation}
\Psi_{k,j,t}(x)=\frac{v_{j}}{\left\Vert e^{itx}\right\Vert }e^{i\left(  2\pi
k+t\right)  x}+O(\frac{\ln|k|)}{k}).
\end{equation}
This formula is uniform with respect to $t\in Q_{\varepsilon}(n),$
$x\in\lbrack0,1],$ that is, there exist a constant $c_{3}$ such that the term
$O(k^{-1}\ln|k|)$ in (8) satisfies%
\[
\mid O(k^{-1}\ln|k|)\mid<c_{3}\mid k^{-1}\ln|k|\mid,\text{ }\forall t\in
Q_{\varepsilon}(n),\text{ }x\in\lbrack0,1].
\]

$(b)$ If $t\in\mathbb{C}(n)$ then the root functions of $L_{t}$ form a Riesz
basis in $L_{2}^{m}(0,1)$.\ 

$(c)$ The eigenfunction $X_{k,j,t}(x)$ of $L_{t}^{\ast},$ where $(X_{k,j,t}%
,\Psi_{k,j,t})=1$, satisfies the following, uniform with respect to $t\in
Q_{\varepsilon}(n),$ $x\in\lbrack0,1],$ formula
\begin{equation}
X_{k,j,t}(x)=v_{j}^{\ast}\left\Vert e^{itx}\right\Vert e^{i(2k\pi+\bar{t}%
)x}+O(\frac{\ln|k|}{k}),
\end{equation}
where $v_{j}^{\ast}$ is the eigenvector of \ $C^{\ast}$ corresponding to
$\overline{\mu_{j}}$ and $(v_{j}^{\ast},v_{j})=1.$

$(d)$ If $f$ is absolutely continuous function satisfying (1) and
$f^{^{\prime}}\in L_{2}^{m}[0,1]$ then the expansion series of $f(x)$ by the
root functions of $L_{t}$ converges uniformly in $[0,1],$ where $t\in
\mathbb{C}(n).$
\end{theorem}

Note that A. A. Shkalikov [13, 14] proved that the root functions of the
operators generated by a ordinary differential expression, in the scalar case,
with summable coefficients and more complicated boundary conditions form a
Riesz basis with brackets. L. M. Luzhina [8] generalized these results for the
matrix case. In [22] we prove that if $n=2$ and the eigenvalues of the matrix
$C$ are simple then the root functions of $L_{t}$ for $t\in(0,\pi)\cup
(\pi,2\pi)$ form a ordinary Riesz basis without brackets. The case $n>2$ is
more complicated and the most part of the method of the paper [22] does not
work here, since in the case $n>2$ the adjoint operator of the operator
generated by $l(y)$ with arbitrary summable coefficients can not be defined by
the Lagrange's formula.

In section 3 using Theorem 2 we obtain spectral expansion for the operator
$L$. The spectral expansion for the Hill operator with real-valued potential
$q(x)$\ was constructed by Gelfand in [4] and Titchmarsh in [15]. Tkachenko
proved in [16] that the Hill operator, namely the operator $L$ in the case
$m=1,$ $n=2$ can be reduced to triangular form if all eigenvalues \ of the
corresponding operators $L_{t}$ for $t\in\lbrack0,2\pi)$ are simple. McGarvey
in [10,11] proved that $L,$ in the case $m=1,$ is spectral operator if the
projections of the operator $L$ are uniformly bounded. Gesztesy and Tkachenko
in the recent paper [5] proved that the Hill operator is a spectral operator
of scalar type if and only if for all $t\in\lbrack0,2\pi)$ the operators
$L_{t}$ have not associated function, the multiple point of either the
periodic or anti-periodic spectrum is a point of its Dirichlet spectrum and
some other condition hold. However, in general, the eigenvalues are not
simple, projections are not uniformly bounded, and $L_{t}$ has associated
function, since the Hill operator with simple potential $q(x)=e^{i2\pi x}$ has
infinitely many spectral singularities ( see [3], where Gasymov investigated
the Hill operator with special potential, analytically continuable onto the
upper half plane). Note that the spectral singularity of $L$ is the point of
$S(T)$ in neighborhood on which the projections of the operator $L$ are not
uniformly bounded and we proved in [18] that a number $\lambda\in
S(L_{t})\subset S(L)$ is a spectral singularity if and only if $L_{t}$ has an
associated function corresponding to the eigenvalue $\lambda.$ The existence
of the spectral singularities and the absence of the Parseval's equality for
the nonself-adjoint operator $L_{t}$ do not allow us to apply the elegant
method of Gelfand ( see\ [4]) for construction of the spectral expansion for
the nonself-adjoin operator $L$. These situation essentially complicate the
construction of\ the spectral expansion for the nonself-adjoint case. In [17]
and [20] we constructed the spectral expansion for the Hill operator with
continuous complex-valued potential $q(x)$ and with locally summable
complex-valued potential $q(x)$ respectively. Then in [19] and [21] we
constructed the spectral expansion for the nonself-adjoint operator $L,$ in
the case $m=1,$ with coefficients $p_{k}\in C^{(k-1)}[0,1]$ \ and with
$p_{k}\in L_{1}[0,1]$ for $k=2,3,...,n$ respectively. In the paper [9] we
constructed the spectral expansion of $L$ when $p_{k,i,j}\in C^{(k-1)}[0,1].$
In this paper we do it when $p_{k,i,j}(x)$ are arbitrary Lebesgue integrable
on $(0,1)$ functions. Besides, in [9] the expansion is obtained for compactly
supported continuous vector functions, while in this paper for each function
$f\in L_{2}^{m}(-\infty,\infty)$ satisfying
\begin{equation}
\sum_{k=-\infty}^{\infty}\mid f(x+k)\mid<\infty
\end{equation}
when $n=2\mu-1$ and for each function from $S,$ where $f(x)\in S\subset
L_{2}^{m}(-\infty,\infty)$ if and only if there exist positive constants $M$
and $\alpha$ such that
\begin{equation}
\mid f(x)\mid<Me^{-\alpha\mid x\mid},\text{ }\forall x\in(-\infty,\infty),
\end{equation}
when $n=2\mu.$ Moreover, using Theorem 2, we prove that the spectral expansion
of $L$ converges uniformly in every bounded subset of $(-\infty,\infty)$ if
$f$ is absolutely continuous compactly supported function and $f^{^{\prime}%
}\in L_{2}^{m}(-\infty,\infty).$ Note that the spectral expansion obtained in
[9], when $p_{k,i,j}\in C^{(k-1)}[0,1],$ converges in the norm of $L_{2}%
^{m}(a,b),$ where $a$ and $b$ are arbitrary real number. Some parts of the
proofs of the spectral expansions for $L$ is just writing in vector form of
the corresponding proofs obtained in [19] for the case $m=1.$ These parts are
given in appendices, in order to give a possibility to reed this paper independently.

\section{On the eigenvalues and root functions of $L_{t}$}

The formula (4) shows that the eigenvalue $\lambda_{k,j}(t)$ of $L_{t}$ is
close to the eigenvalue $\left(  2k\pi i+ti\right)  ^{n}$ of $L_{t}(0).$ If
$t\in Q_{\varepsilon}(n),$ $\mid k\mid\gg1$ then the eigenvalue $\left(  2\pi
ki+ti\right)  ^{n}$ of $L_{t}(0)$ lies far from the other eigenvalues $\left(
2p\pi i+ti\right)  ^{n}.$ It follows from (4) that
\[
|\lambda_{k,j}(t)-\left(  2\pi pi+ti\right)  ^{n}|>c_{4}%
((||k|-|p||+1)(|k|+|p|)^{n-1}%
\]
for $p\neq k,$ $t\in Q_{\varepsilon}(n),$ where $\mid k\mid\gg1$. Using this
one can easily verify that
\begin{equation}
\sum_{p:p>d}\frac{|p|^{n-\nu}}{\left\vert \lambda_{k,j}(t)-(2\pi
pi+ti)^{n}\right\vert }=O(\frac{1}{d^{\nu-1}}),\text{ }\forall d>2\mid
k\mid,\text{ }%
\end{equation}%
\begin{equation}
\sum_{p:p\neq k}\frac{\mid p\mid^{n-\nu}}{\left\vert \lambda_{k,j}(t)-(2\pi
pi+ti)^{n}\right\vert }=O(\frac{\ln|k|}{k^{\nu-1}}),
\end{equation}
where $\mid k\mid\gg1,$ $\nu\geq2,$ and (12), (13) are uniform with respect to
$t\in Q_{\varepsilon}(n)$.

The boundary conditions adjoint to (2) is $U_{\mathbb{\nu}\text{,}\overline
{t}}(y)=0$ for $\nu=0,1,...,(n-1).$ Therefore the eigenfunction $\varphi
_{k,s,t}^{\ast}(x)$ and$\ \Phi_{k,s,t}^{\ast}(x)$ of the operators
$L_{t}^{\ast}(0)$ and $L_{t}^{\ast}(C)$ corresponding to the eigenvalues
$\overline{\left(  2\pi pi+ti\right)  ^{n}}$ and \ $\overline{\mu_{k,j}(t)}$
respectively and satisfying $(\varphi_{k,j,t},\varphi_{k,s,t}^{\ast})=1,$
$(\Phi_{k,j,t},\Phi_{k,s,t}^{\ast})=1$ are
\begin{equation}
\varphi_{k,s,t}^{\ast}(x)=e_{s}\left\Vert e^{itx}\right\Vert e^{i\left(  2\pi
k+\overline{t}\right)  x},\text{ }\Phi_{k,s,t}^{\ast}(x)=v_{s}^{\ast
}\left\Vert e^{itx}\right\Vert e^{i\left(  2\pi k+\overline{t}\right)  x},
\end{equation}
where $v_{s}^{\ast}$ is defined in Theorem 2$(c)$.

To prove the asymptotic formulas for the eigenvalues $\lambda_{k,j}(t)$ and
the corresponding normalized eigenfunctions $\Psi_{k,j,t}(x)$ of $L_{t}$ we
use the formula
\begin{equation}
(\lambda_{k,j}-\mu_{k,s})(\Psi_{k,j,t},\Phi_{k,s,t}^{\ast})=((P_{2}%
-C)\Psi_{k,j,t}^{(n-2)},\Phi_{k,s,t}^{\ast})+\sum_{\nu=3}^{n}(P_{\nu}%
\Psi_{k,j,t}^{(n-\nu)},\Phi_{k,s,t}^{\ast})
\end{equation}
which can be obtained from
\begin{equation}
L_{t}\Psi_{k,j,t}(x)=\lambda_{k,j}(t)\Psi_{k,j,t}(x)
\end{equation}
by multiplying scalarly by $\Phi_{k,s,t}^{\ast}(x)$. To estimate the
right-hand side of (15)\ we use (12), (13), the following lemma, and the
formula
\begin{equation}
\left(  \lambda_{k,j}(t)-\left(  2\pi pi+ti\right)  ^{n}\right)  \left(
\Psi_{k,j,t},\varphi_{p,s,t}^{\ast}\right)  =\sum_{\nu=2}^{n}(P_{\nu}%
\Psi_{k,j,t}^{(n-\nu)},\varphi_{p,s,t}^{\ast})
\end{equation}
which can be obtained from (16) by multiplying scalarly by $\varphi
_{p,s,t}^{\ast}(x)$.

\begin{lemma}
If $|k|\gg1$ and $t\in Q_{\varepsilon}(n),$ then
\begin{equation}
(P_{\nu}\Psi_{k,j,t}^{(n-\nu)},\varphi_{p,s,t}^{\ast})=\sum_{q=1}^{m}%
(\sum_{l=-\infty}^{\infty}p_{\nu,s,q,p-l}(2\pi li+it)^{n-\nu}(\Psi
_{k,t},\varphi_{l,q,t}^{\ast})),
\end{equation}
where $p_{\nu,s,q,k}=\int_{0}^{1}p_{\nu,s,q}(x)e^{-i2\pi kx}dx.$ Moreover
\begin{equation}
\max_{p\in\mathbb{Z}\text{,}s=1,2,...,m}\left\vert \sum_{\nu=2}^{n}(P_{\nu
}\Psi_{k,j,t}^{(n-\nu)},\varphi_{p,s,t}^{\ast})\right\vert <c_{5}|k|^{n-2}.
\end{equation}

\end{lemma}

\begin{proof}
Since $P_{2}\Psi_{k,j,t}^{(n-2)}+P_{3}\Psi_{k,j,t}^{(n-3)}+\cdots+P_{n}%
\Psi_{k,j,t}\in L_{1}^{m}[0,1]$ we have
\[
\lim_{p\rightarrow\infty}\left\vert \sum_{\nu=2}^{n}(P_{\nu}\Psi
_{k,j,t}^{(n-\nu)},\varphi_{p,s,t}^{\ast})\right\vert =0.
\]
Therefore there exist a positive constant $M(k,j)$ and indices $p_{0},s_{0}$
satisfying
\begin{equation}
\max_{\substack{p\in\mathbb{Z}\text{,}\\s=1,2,...,m}}\left\vert \sum_{\nu
=2}^{n}(P_{\nu}\Psi_{k,j,t}^{(n-\nu)},\varphi_{p,s,t}^{\ast})\right\vert
=\left\vert \sum_{\nu=2}^{n}(P_{\nu}\Psi_{k,j,t}^{(n-\nu)},\varphi
_{p_{0},s_{0},t}^{\ast})\right\vert =M(k,j).
\end{equation}
Then using (17) and (12), we get
\begin{equation}
\left\vert \left(  \Psi_{k,j,t},\varphi_{p,s,t}^{\ast}\right)  \right\vert
\leq\frac{M(k,j)}{\left\vert \lambda_{k,j}(t)-\left(  2\pi pi+it\right)
^{n}\right\vert },
\end{equation}%
\[
\sum_{p:|p|>d}\left\vert \left(  \Psi_{k,j,t},\varphi_{p,s,t}^{\ast}\right)
\right\vert <\frac{c_{6}M(k,j)}{d^{n-1}},
\]
where $d>2|k|$. This implies that the decomposition of $\Psi_{k,j,t}(x)$ by basis

$\{\varphi_{p,s,t}(x):p\in\mathbb{Z}$, $s=1,2,...,m\}$ is of the form
\begin{equation}
\Psi_{k,j,t}(x)=\sum_{p:|p|\leq d}\left(  \Psi_{k,j,t},\varphi_{p,s,t}^{\ast
}\right)  \varphi_{p,s,t}(x)+g_{0,d}(x),
\end{equation}
where%
\[
\sup_{x\in\lbrack0,1]}|g_{0,d}(x)|<\frac{c_{7}M(k,j)}{d^{n-1}}.
\]
Now using the integration by parts, (1), and the inequality (21), we obtain%
\[
(\Psi_{k,j,t}^{(n-\nu)},\varphi_{p,s,t}^{\ast})=(2\pi ip+it)^{n-\nu}\left(
\Psi_{k,j,t},\varphi_{p,s,t}^{\ast}\right)  ,
\]%
\[
\left\vert \left(  \Psi_{k,j,t}^{(n-\nu)},\varphi_{p,s,t}^{\ast}\right)
\right\vert \leq\frac{|2\pi ip+it|^{n-\nu}M(k,j)}{\left\vert \lambda
_{k}(t)-(2\pi pi+it)^{n}\right\vert }.
\]
Therefore arguing as in the proof of (22) and using (12) we get%
\begin{equation}
\Psi_{k,j,t}^{(n-\nu)}(x)=\sum_{p:|p|\leq d}\left(  \Psi_{k,j,t}^{(n-\nu
)},\varphi_{p,s,t}^{\ast}\right)  \varphi_{p,s,t}(x)+g_{\nu,d}(x),
\end{equation}
where $\nu=2,3,\ldots,n$, and%
\[
\sup_{x\in\lbrack0,1]}|g_{\nu,d}(x)|<\frac{c_{7}M(k,j)}{d^{\nu-1}}.
\]
Now using (23) in $(P_{\nu}\Psi_{k,j,t}^{(n-\nu)},\varphi_{p,s,t}^{\ast})$ and
tending $q$ to $\infty$, we obtain (18).

Let us we prove (19). It follows from (20) and (18) that
\begin{align}
&  M(k,j)=\left\vert \sum_{\nu=2}^{n}(P_{\nu}\Psi_{k,j,t}^{(n-\nu)}%
,\varphi_{p_{0},s_{0},t}^{\ast})\right\vert \nonumber\label{f2}\\
&  =\left\vert \sum_{\nu=2}^{n}\sum_{q=1}^{m}(\sum_{l=-\infty}^{\infty}%
p_{\nu,s_{0},q,p_{0}-l}(2\pi im+it)^{n-\nu}(\Psi_{k,j,t},\varphi_{l,q,t}%
^{\ast}))\right\vert .
\end{align}
By (21) and (13) we have
\[
\left\vert \sum_{\nu=2}^{n}\sum_{q=1}^{m}(\sum_{l\neq k}p_{\nu,s_{0}%
,q,p_{0}-l}(2\pi im+it)^{n-\nu}(\Psi_{k,j,t},\varphi_{l,q,t}^{\ast
}))\right\vert \leq c_{8}M(k,j)\frac{\ln|k|}{\mid k\mid}.
\]
On the other hand
\[
\left\vert \sum_{\nu=2}^{n}\sum_{q=1}^{m}(p_{\nu,s_{0},q,p_{0}-k}(2\pi
im+it)^{n-\nu}(\Psi_{k,j,t},\varphi_{k,q,t}^{\ast}))\right\vert =O(k^{n-2}).
\]
Therefore using (24) we get
\[
M(k,j)=M(k,j)O((\frac{\ln|k|}{k})+O(|k|^{n-2}),
\]
$M(k,j)=O(|k|^{n-2})$ which means that (19) holds
\end{proof}

It follows from (19)-(21) that
\begin{equation}
\left\vert \left(  \Psi_{k,j,t},\varphi_{p,q,t}^{\ast}\right)  \right\vert
\leq\frac{c_{5}|k|^{n-2}}{\left\vert \lambda_{k,j}(t)-(2\pi pi+it)^{n}%
\right\vert },\text{ }\forall p\neq k.
\end{equation}

Now using this we prove the following lemma.

\begin{lemma}
The following equalities
\begin{equation}
\left(  (P_{2}-C)\Psi_{k,j,t}^{(n-2)},\Phi_{k,s,t}^{\ast}\right)
=O(k^{n-3}\ln|k|),
\end{equation}%
\begin{equation}
\left(  (P_{\nu}\Psi_{k,j,t}^{(n-\nu)},\Phi_{k,s,t}^{\ast}\right)  =O(k^{n-3})
\end{equation}
hold uniformly with respect to $t\in Q_{\varepsilon}(n),$ where $\nu\geq3.$
\end{lemma}

\begin{proof}
Using (18) for $\nu=2$, $p=k$ and the obvious relation
\begin{equation}
\left(  C\Psi_{k,j,t}^{(n-2)},\varphi_{k,s,t}^{\ast}\right)  =\sum_{q=1}%
^{m}(p_{2,s,q,0}(2\pi ki+it)^{n-2}(\Psi_{k,j,t},\varphi_{k,q,t}^{\ast
}))\nonumber
\end{equation}
we see that
\[
\left(  (P_{2}-C)\Psi_{k,j,t}^{(n-2)},\varphi_{k,s,t}^{\ast}\right)
=\sum_{q=1}^{m}(\sum_{l\neq k}p_{2,s,q,k-l}(2\pi li+it)^{n-2}(\Psi
_{k,j,t},\varphi_{l,q,t}^{\ast})).
\]
This with (25) and (13) for $\nu=2$ implies that
\[
\left(  (P_{2}-C)\Psi_{k,j,t}^{(n-2)},\varphi_{k,s,t}^{\ast}\right)
=O(k^{n-3}\ln|k|).
\]
Similarly, using (18), (25), (13) we obtain
\[
\left(  (P_{\nu}\Psi_{k,j,t}^{(n-\nu)},\varphi_{k,s,t}^{\ast}\right)
=O(k^{n-3}),\text{ }\forall\nu\geq3.
\]
Since (13) is uniform with respect to $t\in Q_{\varepsilon}(n)$ and the
constant $c_{5}$ in (25) does not depend on $t$ ( recall that we denote by
$c_{k}$ the constant independent of $t$) these formulas are uniform with
respect to $t\in Q_{\varepsilon}(n).$ Therefore recalling the definitions of
$\Phi_{k,s,t}^{\ast}$ and $\varphi_{k,q,t}^{\ast}$ ( see (14)) we get the
proof of (26) and (27)
\end{proof}

\begin{lemma}
There exist positive number $N_{1},$ independent of $t,$ such that%
\begin{equation}
\max_{s=1,2,...,m}\left\vert \left(  \Psi_{k,j,t},\Phi_{k,s,t}^{\ast}\right)
\right\vert >c_{9}%
\end{equation}
for all $\mid k\mid\geq N_{1},$ $t\in Q_{\varepsilon}(n),$ and $j=1,2,...,m.$
\end{lemma}

\begin{proof}
It follows from (25) and (13) that
\begin{equation}
\sum\limits_{s=1,2,...,m}(\sum_{p:\text{ }p\neq k}\left\vert \left(
\Psi_{k,j,t},\varphi_{p,s,t}^{\ast}\right)  \right\vert =O(\frac{\ln|k|}{k})
\end{equation}
and this formula is uniform with respect to $t\in Q_{\varepsilon}(n).$ Then
the decomposition of $\Psi_{k,j,t}(x)$ by the basis $\{\varphi_{p,s,t}%
(x):s=1,2,...,m$, $p\in\mathbb{Z\}}$ has the form%
\begin{equation}
\Psi_{k,j,t}(x)=\sum\limits_{s=1,2,...,m}\left(  \Psi_{k,j,t},\varphi
_{k,s,t}^{\ast}\right)  \varphi_{k,s,t}(x)+O(\frac{\ln|k|}{k}).
\end{equation}
Since $\parallel\Psi_{k,j,t}\parallel=\parallel\varphi_{k,j,t}\parallel=1$ and
(30) is uniform with respect to $t\in Q_{\varepsilon}(n),$ there exists a
positive constant $N_{1},$ independent of $t,$ such that%
\[
\max_{s=1,2,...,m}\left\vert \left(  \Psi_{k,j,t},\varphi_{k,s,t}^{\ast
}\right)  \right\vert >\frac{1}{m+1}%
\]
for all $\mid k\mid\geq N_{1},$ $t\in Q_{\varepsilon}(n)$ and $j=1,2,...,m.$
Therefore using (14) and taking into account that the vectors $v_{1}^{\ast
},v_{2}^{\ast},...,v_{m}^{\ast}$ form a basis in $\mathbb{C}^{m},$ that is,
$e_{s}$ is a linear combination of these vectors we get the proof of (28)
\end{proof}

\textbf{THE PROOF OF THEOREM 1}$(a)$\textbf{.} It follows from Lemma 2 that
there exists a positive constant $N_{2},$ independent of $t,$ such that if
$\mid k\mid\geq N_{2},$ $t\in Q_{\varepsilon}(n)$ then the right-hand side of
(15) is less than $c_{10}|k|^{n-3}\ln|k|.$ Therefore (15) and Lemma 3 give the
proof of the Theorem 1$(a)$.

\textbf{THE PROOF OF THEOREM 1}$(b)$\textbf{. }Let $\lambda_{k,j}$ be an
eigenvalue of $L_{t}$ lying in $U(\mu_{k,p(j)}(t),c_{1}|k|^{n-3}\ln|k|)$ and
$\Psi_{k,j,t}$ be any normalized eigenfunction corresponding to $\lambda
_{k,j}.$ Then using (5) and taking into account that the eigenvalues of $C$
are simple we get
\[
\mid\lambda_{k,j}-\mu_{k,s}\mid>a_{p(j)}\mid k\mid^{n-2}\text{ \ for }s\neq
p(j),
\]
where $a_{p(j)}=\min_{s\neq p(j)}\mid\mu_{p(j)}-\mu_{s}\mid.$ This with (15),
(26), (27) gives
\begin{equation}
\left(  \Psi_{k,j,t},\Phi_{k,s,t}^{\ast}\right)  =O(\frac{(\ln|k|)}{k}),\text{
}\forall s\neq p(j).
\end{equation}
On the other hand by (14) and (29) we have
\begin{equation}
\sum\limits_{s=1,2,...,m}(\sum_{p:\text{ }p\neq k}\left\vert \left(
\Psi_{k,j,t},\Phi_{p,s,t}^{\ast}\right)  \right\vert =O(\frac{(\ln|k|)}{k}).
\end{equation}
Since (26), (27), (29) are uniform with respect to $t\in Q_{\varepsilon}(n)$
the formulas (31) and (32) are also uniform. Therefore decomposing
$\Psi_{k,j,t}(x)$ by basis $\{\Phi_{p,s,t}(x):s=1,2,...,m$, $p\in\mathbb{Z\}}$
we see that any normalized eigenfunction corresponding to $\lambda_{k,j}$
satisfies (6). If there are two linearly independent eigenfunctions
corresponding to $\lambda_{k,j},$ then one can find two orthogonal
eigenfunctions satisfying (6), which is impossible. Theorem 1 is proved.

To proof of the main results for $L_{t}$ (Theorem 2) we need to investigate
the normalized associated function $\Psi_{k,j,1,t}(x)$ of $L_{t}$
corresponding to the eigenvalue $\lambda_{k,j}(t)$. By definition of the
associated function we have
\begin{equation}
(L_{t}-\lambda_{k,j})\Psi_{k,j,1,t}(x)=\Psi_{k,j,0,t}(x),
\end{equation}
where $\Psi_{k,j,0,t}(x)$ is an eigenfunction of $L_{t}.$ Note that, in
general, the eigenfunction $\Psi_{k,j,0,t}(x)$ is not normalized. For
investigation of the associated function we use the following formulas.
Multiplying scalarly (33) by $\varphi_{p,s,t}^{\ast}$ we get
\begin{equation}
(\lambda_{k,j}-\left(  2\pi pi+ti\right)  ^{n})(\Psi_{k,j,s,t},\varphi
_{p,s,t}^{\ast})=\sum_{\nu=2}^{n}(P_{\nu}\Psi_{k,j,q,t}^{(n-\nu)}%
,\varphi_{p,s,t}^{\ast})-(\Psi_{k,j,0,t},\varphi_{p,s,t}^{\ast}).
\end{equation}
Similarly, multiplying scalarly (33) by $\Phi_{k,s,t}^{\ast},$ we obtain
\[
(\Psi_{k,j,0,t},\Phi_{k,s,t}^{\ast})=(L_{t}(C)\Psi_{k,j,1,t},\Phi
_{k,s,t}^{\ast})+((P_{2}-C)\Psi_{k,j,1,t}^{(n-2)},\Phi_{k,s,t}^{\ast})+
\]%
\[
\sum_{\nu=3}^{n}(P_{\nu}\Psi_{k,j,1,t}^{(n-\nu)},\Phi_{k,s,t}^{\ast}%
)-\lambda_{k,j}(\Psi_{k,j,1,t},\Phi_{k,s,t}^{\ast}).
\]
Since $(L_{t}(C)\Psi_{k,j,1,t},\Phi_{k,s,t}^{\ast})=\mu_{k,s}(\Psi
_{k,j,1,t},\Phi_{k,s,t}^{\ast})$ we have
\[
(\lambda_{k,j}-\mu_{k,s})(\Psi_{k,j,1,t},\Phi_{k,s,t}^{\ast})=
\]%
\begin{equation}
((P_{2}-C)\Psi_{k,j,1,t}^{(n-2)},\Phi_{k,s,t}^{\ast})+\sum_{\nu=3}^{n}(P_{\nu
}\Psi_{k,j,1,t}^{(n-\nu)},\Phi_{k,s,t}^{\ast})-(\Psi_{k,j,0,t},\Phi
_{k,s,t}^{\ast}).
\end{equation}

\begin{lemma}
For any normalized associated eigenfunction $\Psi_{k,j,1,t}$ of $L_{t}$ the
following, uniform with respect to $t\in Q_{\varepsilon}(n),$ formulas hold
\begin{equation}
\left(  (P_{2}-C)\Psi_{k,j,1,t}^{(n-2)},\Phi_{k,s,t}^{\ast}\right)
=O(k^{n-3}\ln|k|),
\end{equation}%
\begin{equation}
\left(  (P_{\nu}\Psi_{k,j,1,t}^{(n-\nu)},\Phi_{k,s,t}^{\ast}\right)
=O(k^{n-3}),\text{ }\forall\nu\geq3.
\end{equation}

\end{lemma}

\begin{proof}
Instead of (17) using (34) and repeating the proof of (19) we obtain
\begin{equation}
\max_{p\in\mathbb{Z}\text{,}s=1,2,...,m}\left\vert \sum_{\nu=2}^{n}(P_{\nu
}\Psi_{k,j,q,t}^{(n-\nu)},\varphi_{p,s,t}^{\ast})\right\vert <c_{11}%
(|k|^{n-2}+\parallel\Psi_{k,j,0,t}\parallel).
\end{equation}
Using (38) and repeating the proof of (25)-(27) we get
\begin{equation}
\left\vert (\Psi_{k,j,1t},\varphi_{p,s,t}^{\ast})\right\vert \leq\frac
{c_{12}(|k|^{n-2}+\parallel\Psi_{k,j,0,t}\parallel)}{\left\vert \lambda
_{k}(t)-(2\pi pi+it)^{n}\right\vert },
\end{equation}%
\begin{equation}
((P_{2}-C)\Psi_{k,j,1,t}^{(n-2)},\Phi_{k,s,t}^{\ast})=O(|k|^{n-3}%
\ln|k|+\parallel\Psi_{k,j,0,t}\parallel\frac{\ln|k|}{|k|})),
\end{equation}%
\begin{equation}
\left(  (P_{\nu}\Psi_{k,j,1,t}^{(n-\nu)},\Phi_{k,s,t}^{\ast}\right)
=O(|k|^{n-3}+|k|^{-1}\parallel\Psi_{k,j,0,t}\parallel)
\end{equation}
for $\nu\geq3.$ Using (40), (41) in (35) for $\ s=p(j)$ and taking into
account that
\begin{align*}
(\lambda_{k,j}-\mu_{k,p(j)})(\Psi_{k,j,1,t},\Phi_{k,p(j),t}^{\ast})  &
=O(\frac{\ln|k|}{|k|^{3-n}}),\\
(\frac{\Psi_{k,j,0,t}}{\parallel\Psi_{k,j,0,t}\parallel},\Phi_{k,p(j),t}%
^{\ast})  &  =1+O(\frac{\ln|k|}{k})
\end{align*}
(see the definition of $p(j)$ in Theorem 1) we obtain
\[
O(\frac{\ln|k|}{|k|^{3-n}})=\parallel\Psi_{k,j,0,t}\parallel(1+O(\frac{\ln
|k|}{k}))+O(\frac{\ln|k|}{|k|^{3-n}}+\parallel\Psi_{k,j,0,t}\parallel\frac
{\ln|k|}{|k|})
\]
which yields the equality
\begin{equation}
\parallel\Psi_{k,j,0,t}(x)\parallel=O(|k|^{n-3}\ln|k|).
\end{equation}
Now (40), (41) and (42) imply (36) and (37)
\end{proof}

\begin{lemma}
Any normalized associated function $\Psi_{k,j,1,t}(x)$ of $L_{t}$
corresponding to the eigenvalue $\lambda_{k,j}(t)\in U(\mu_{k,p(j)}%
(t),c_{1}|k|^{n-3}\ln|k|),$ where$\mid k\mid\geq N_{0}$ and $c_{1},$ $p(j),$
$N_{0}$ are defined in Theorem 1, satisfies
\begin{equation}
\Psi_{k,j,1,t}(x)=\frac{v_{p(j)}}{\left\Vert e^{itx}\right\Vert }e^{i\left(
2\pi k+t\right)  x}+O(\frac{(\ln|k|)}{k}).
\end{equation}

\end{lemma}

\begin{proof}
$(a)$ It follows from (39), (42) that
\[
\left\vert \left(  \Psi_{k,j,1t},\varphi_{p,s,t}^{\ast}\right)  \right\vert
\leq\frac{c_{13}|k|^{n-2}}{\left\vert \lambda_{k}(t)-(2\pi pi+it)^{n}%
\right\vert }.
\]
Using this instead of (25) and repeating the proof of (32) and (31) we obtain
\[
\sum\limits_{s=1,2,...,m}(\sum_{p:p\neq k}\left\vert \left(  \Psi
_{k,j,1,t},\Phi_{p,s,t}^{\ast}\right)  \right\vert =O(\frac{\ln|k|}{k}),
\]%
\[
\left(  \Psi_{k,j,1,t},\Phi_{k,s,t}^{\ast}\right)  =O(\frac{\ln|k|}{k})\text{,
}\forall s\neq p(j)
\]
which imply the proof of (43)
\end{proof}

\textbf{THE PROOF OF THEOREM 2}$(a).$ Let $\lambda_{k,j}(t)$ be eigenvalue of
$L_{t}$ lying in

$U(\mu_{k,p(j)}(t),c_{1}|k|^{n-3}\ln|k|),$ where $\mid k\mid\geq N_{0}.$ By
Theorem 1 there exist only one eigenfunction $\Psi_{k,j,t}(x)$ corresponding
to $\lambda_{k,j}(t).$ Suppose that there exist associated function
$\Psi_{k,j,1,t}(x)$ corresponding to the eigenvalue $\lambda_{k,j}(t).$ Using
Lemma 5 and taking into account that for any $a\in\mathbb{C}$ the function
$\Psi_{k,j,1t}+a\Psi_{k,j,t}$ is associated function one can find two
orthogonal root functions satisfying (43) which is impossible. Thus we proved
that the operator $L_{t}$ has not associated function corresponding to the
eigenvalue $\lambda_{k,j}(t)$ for $\mid k\mid\geq N_{0}.$ Using this, (3),
(4), and Theorem 1, we obtain the following:

\begin{proposition}
There exist a number $N_{0}$ such that the eigenvalues $\lambda_{k,1}(t),$
$\lambda_{k,2}(t),$ ..., $\lambda_{k,m}(t)$ of $L_{t}$ for $t\in
Q_{\varepsilon}(n),$ $\mid k\mid\geq N_{0}$ are simple and they lie in the
union of the pairwise disjoint intervals
\begin{equation}
U(\mu_{k,1}(t),\frac{c_{1}\ln|k|}{|k|^{3-n}}),\text{ }U(\mu_{k,2}%
(t),\frac{c_{1}\ln|k|}{|k|^{3-n}}),...,\text{ }U(\mu_{k,m}(t),\frac{c_{1}%
\ln|k|}{|k|^{3-n}}).
\end{equation}

\end{proposition}

Now let us prove that in each of these intervals there exists unique
eigenvalue of $L_{t}.$ For this we consider the following family of operators
\begin{equation}
L_{t,\varepsilon}=L_{t}(C)+\varepsilon(L_{t}-L_{t}(C)),\text{ }0\leq
\varepsilon\leq1.
\end{equation}
It is clear that the proposition 1 holds for $L_{t,\varepsilon},$ that is, the
eigenvalues $\lambda_{k,1,\varepsilon}(t),$ $\lambda_{k,2,\varepsilon}(t),$
..., $\lambda_{k,m,\varepsilon}(t),$ where $|k|\geq N_{0},$ of
$L_{t,\varepsilon}$ are simple and they lie in the union of the pairwise
disjoint $m$ intervals (44). Since $\lambda_{k,j,\varepsilon}$ is a simple
eigenvalue it is a simple root of the characteristic determinant
$\Delta(\lambda,\varepsilon)$ of $L_{t,\varepsilon}.$ Clearly, $\Delta
(\lambda,\varepsilon)$ is analytic function of $\lambda$ and $\varepsilon$ and
$\Delta(\lambda_{k,j,\varepsilon},\varepsilon)=0,$ $\frac{\partial}%
{\partial\lambda}\Delta(\lambda,\varepsilon)\neq0$ for $\lambda=\lambda
_{k,j,\varepsilon}.$ Therefore using the implicit function theorem and taking
into account that $\lambda_{k,1,\varepsilon}(t),$ $\lambda_{k,2,\varepsilon
}(t),$ ..., $\lambda_{k,m,\varepsilon}(t)$ are simple eigenvalues one can
easily see that these eigenvalues continuously depend on $\varepsilon.$
Therefore taking into account that in each of the pairwise disjoint intervals
(44) there exists unique eigenvalue of $L_{t,0},$ we conclude that in
$U(\mu_{k,j},\frac{c_{1}\ln|k|}{|k|^{3-n}})$ for $|k|\geq N_{0},$
$j=1,2,...,m$ there exists unique eigenvalue of $L_{t,\varepsilon}$ for all
values of $\varepsilon\in\lbrack0,1].$ Let us denote this eigenvalue of
$L_{t,\varepsilon}$ by $\lambda_{k,j,\varepsilon}(t).$ Thus we proved the following:

\begin{proposition}
Let $t\in Q_{\varepsilon}(n),$ $\varepsilon\in\lbrack0,1].$ All large
eigenvalues of $L_{t,\varepsilon}$ belong to one of the intervals (44) for
$|k|\geq N_{0}.$ For each eigenvalues $\mu_{k,j}(t)$ of $L_{t}\left(
C\right)  ,$ where $|k|\geq N_{0},$ there exists unique eigenvalue
$\lambda_{k,j,\varepsilon}(t)$ of $L_{t,\varepsilon}$ lying in $U(\mu
_{k,j}(t),c_{1}|k|^{n-3}\ln|k|).$
\end{proposition}

By Proposition 1 the eigenvalue $\lambda_{k,j}(t)$ of $L_{t}$ for $|k|\geq
N_{0}$ is simple and by Theorem 1 the corresponding eigenfunction satisfy (6),
where $p(j)=j$ (see the definition of $p(j)$ in Theorem 1), that is, (8), (7)
and Theorem 2$(a)$ is proved.

\textbf{THE PROOF OF THEOREM 2}$(b).$ It follows from (8) that the root
functions of $L_{t}$ quadratically close to the system
\[
\{v_{j}\left\Vert e^{itx}\right\Vert ^{-1}e^{i\left(  2\pi k+t\right)
x}:\text{ }k\in\mathbb{Z}\text{, }l=1,2,...,m\}
\]
which form a Riesz in $L_{2}^{m}(0,1).$ On the other hand the system of the
root functions of $L_{t}$ is complete and minimal in $L_{2}^{m}(0,1)$ ( see
[8]). Therefore, by Bari theorem ( see [1,6]), the system of the root
functions of $L_{t}$ forms a Riesz basis in $L_{2}^{m}(0,1).$

\textbf{THE PROOF OF THEOREM 2}$(c).$ To prove the asymptotic formulas for
normalized eigenfunction $\Psi_{k,j,t}^{\ast}(x)$ of $L_{t}^{\ast}$
corresponding to the eigenvalue $\overline{\lambda_{k,j}(t)}$ we use the
formula
\[
\left(  \overline{\lambda_{k,j}(t)}-\overline{\left(  2\pi pi+ti\right)  ^{n}%
}\right)  \left(  \Psi_{k,j,t}^{\ast},\varphi_{p,s,t}\right)  =\sum_{\nu
=2}^{n}(\Psi_{k,j,t}^{\ast},\overline{\left(  2\pi pi+ti\right)  ^{n-\nu}%
}P_{\nu}\varphi_{p,s,t})
\]
obtained from $L_{t}^{\ast}\Psi_{k,j,t}^{\ast}=\overline{\lambda_{k,j}(t)}%
\Psi_{k,j,t}^{\ast}$ by multiplying by $\varphi_{p,s,t}$ and using
\[
(L_{t}^{\ast}\Psi_{k,j,t}^{\ast},\varphi_{p,s,t})=(\Psi_{k,j,t}^{\ast}%
,L_{t}\varphi_{p,s,t}).
\]
Instead of (17) using these formula and arguing as in the proof of (25) we
obtain
\[
\left\vert \left(  \Psi_{k,j,t}^{\ast},\varphi_{p,q,t}\right)  \right\vert
\leq\frac{c_{14}|k|^{n-2}}{\left\vert \lambda_{k,j}(t)-(2\pi pi+it)^{n}%
\right\vert },\text{ }\forall p\neq k.
\]
This with (5) and (13) implies the following relations
\begin{equation}
\left\vert \left(  \Psi_{k,j,t}^{\ast},\Phi_{p,q,t}\right)  \right\vert
\leq\frac{c_{15}|k|^{n-2}}{\left\vert \lambda_{k,j}(t)-(2\pi pi+it)^{n}%
\right\vert },\text{ }\forall p\neq k,
\end{equation}%
\begin{equation}
\sum\limits_{s=1,2,...,m}(\sum_{p:\text{ }p\neq k}\left\vert \left(
\Psi_{k,j,t}^{\ast},\Phi_{p,s,t}\right)  \right\vert =O(\frac{(\ln|k|)}{k}).
\end{equation}
On the other hand (8) and equality $\left(  \Psi_{k,j,t}^{\ast},\Psi
_{k,s,t}\right)  =0$ for $j\neq s$ give
\begin{equation}
\left(  \Psi_{k,j,t}^{\ast},\Phi_{k,s,t}\right)  =O(\frac{(\ln|k|)}{k}),\text{
}\forall s\neq j.
\end{equation}
Since (8), (13) hold uniformly the formulas (46)-(48) are uniform with respect
to $t\in Q_{\varepsilon}(n)$ and they yield
\begin{equation}
\Psi_{k,j,t}^{\ast}(x)=v_{j}^{\ast}\left\Vert e^{itx}\right\Vert e^{(2k\pi
i+i\bar{t})x}+O(\frac{\ln\mid k\mid}{k}),
\end{equation}
where $v_{j}^{\ast}$ is defined in Theorem 2$(c).$ Now (8) and (49) imply (9),
since
\begin{equation}
X_{k,j,t}=\frac{\Psi_{k,j,t}^{\ast}}{(\Psi_{k,j,t}^{\ast},\Psi_{k,j,t}%
)}=(1+O(\frac{\ln\mid k\mid}{k}))\Psi_{k,j,t}^{\ast}.
\end{equation}

\textbf{THE PROOF OF THEOREM 2}$(d).$ To investigate the convergence of the
expansion series of $L_{t}$ we consider the series
\begin{equation}
\sum_{k:\mid k\mid\geq N\text{, }j=1,2,...,m}(f,X_{k,j,t})\Psi_{k,j,t}(x),
\end{equation}
where $N\geq N_{0}$ and $N_{0}$ is defined in Theorem 1, $f(x)$ is absolutely
continuous function satisfying (1) and $f^{^{\prime}}(x)\in L_{2}^{m}(0,1)$.
Without loss of generality instead of the series (51) we consider the series
\begin{equation}
\sum_{k:\mid k\mid\geq N\text{, }j=1,2,...,m}(f_{t},X_{k,j,t})\Psi_{k,j,t}(x),
\end{equation}
where $f_{t}(x)$ is defined by Gelfand transform ( see [4])
\begin{equation}
f_{t}(x)=\sum_{k=-\infty}^{\infty}f(x+k)e^{-ikt},
\end{equation}
$f$ is absolutely continuous compactly supported function and $f^{^{\prime}%
}\in L_{2}^{m}(-\infty,\infty),$ since we use (52) in next section for
spectral expansion of $L.$ It follows from (53) that
\begin{equation}
\text{ }f_{t}(x+1)=e^{it}f_{t}(x),\text{ }f_{t}^{^{\prime}}\in L_{2}^{m}[0,1].
\end{equation}
To prove the uniform convergence of (52) we consider the series%
\begin{equation}
\sum_{\mid k\mid\geq N\text{, }j=1,2,...,m}\mid(f_{t},X_{k,j,t})\mid.
\end{equation}
To estimate the terms of this series we decompose $X_{k,j,t}$ by basis

$\{\Phi_{p,s,t}^{\ast}:p\in\mathbb{Z}$, $s=1,2,...,m\}$ and then use the
inequality
\begin{align}
&  \mid(f_{t},X_{k,j,t})\mid\leq\sum_{\text{ }s=1,2,...,m}\mid(f_{t}%
,\Phi_{k,s,t}^{\ast})\mid\left\vert \left(  X_{k,j,t},\Phi_{k,s,t}\right)
\right\vert +\\
\sum_{p\neq k,\text{ }s=1,2,...,m}  &  \mid(f_{t},\Phi_{p,j,t}^{\ast}%
)\mid\left\vert \left(  X_{k,j,t},\Phi_{p,s,t}\right)  \right\vert .\nonumber
\end{align}
Using the integration by parts and then Schwarz inequality we get%
\begin{equation}
\sum_{\substack{\mid k\mid\geq N\text{,}\\\text{ }s=1,2,...,m}}\mid(f_{t}%
,\Phi_{k,s,t}^{\ast})\mid=\sum_{\substack{\mid k\mid\geq N\text{,}\\\text{
}s=1,2,...,m}}\mid\frac{1}{2\pi ki+it}(f_{t}^{^{\prime}},\Phi_{k,s,t}^{\ast
})\mid<\infty.
\end{equation}
Again using the integration by parts, Schwarz inequality and (46), (50) we
obtain that the expression in the in the second row of (56) is less than
\[
c_{16}\parallel f_{t}^{^{\prime}}\parallel\left(  \sum_{p\neq k,\text{
}s=1,2,...,m}\left\vert \frac{1}{p}\frac{|k|^{n-2}}{\left\vert \lambda
_{k,s}(t)-(2\pi pi+it)^{n}\right\vert }\right\vert ^{2}\right)  ^{\frac{1}{2}%
}.
\]
It is not hard to see that this expression is less than $c_{17}k^{-2},$ that
is, the expression in the second row of (56) is less than $c_{17}k^{-2}.$
Therefore the relations (56), (57) imply that the expressions in (55) and (52)
tend to zero uniformly with respect to $t\in Q_{\varepsilon}(n)$ and $t\in
Q_{\varepsilon}(n),$ $x\in\lbrack0,1]$ respectively as $N\rightarrow\infty.$
Since in the proof of the uniform convergence of (52) we used only the
properties (54) of $f_{t}$ the series (51) converges uniformly with respect to
$x\in\lbrack0,1],$ that is, Theorem 2($d$) is proved.

Note that in the proof of Theorem 2$(d)$ we proved the following theorem,
which will be used in next section.

\begin{theorem}
If $f$ is absolutely continuous, compactly supported function and
$f^{^{\prime}}\in L_{2}^{m}(-\infty,\infty)$ then the series (52), where
$f_{t}$ is defined by (53), $N\geq N_{0}$, $N_{0}$ is defined in Theorem
1$(a),$ converges uniformly with respect to $t\in Q_{\varepsilon}(n),$ $x\in
D$ for any bounded subset $D$ of $(-\infty,\infty).$
\end{theorem}

Indeed we proved that (52) converges uniformly with respect to $t\in
Q_{\varepsilon}(n),$ $x\in\lbrack0,1].$ Therefore taking into account that (1)
implies the equality%
\begin{equation}
\Psi_{k,j,t}(x+1)=e^{it}\Psi_{k,j,t}(x),
\end{equation}
we get the proof of Theorem 3.

\section{Spectral Expansion for $L$}

Let $Y_{1}(x,\lambda),Y_{2}(x,\lambda),\ldots,Y_{n}(x,\lambda)$ be the
solutions of the matrix equation
\begin{equation}
Y^{(n)}(x)+P_{2}\left(  x\right)  Y^{(n-2)}(x)+P_{3}\left(  x\right)
Y^{(n-3)}(x)+...+P_{n}(x)Y=\lambda Y(x),
\end{equation}
satisfying $Y_{k}^{(j)}(0,\lambda)=0_{m}$ for $j\neq k-1$ and $Y_{k}%
^{(k-1)}(0,\lambda)=I_{m},$ where $0_{m}$ and $I_{m}$ are $m\times m$ zero and
identity matrices respectively. The eigenvalues of the operator $L_{t}$ are
the roots of the characteristic determinant
\begin{equation}
\Delta(\lambda,t)=\det(Y_{j}^{(\nu-1)}(1,\lambda)-e^{it}Y_{j}^{(\nu
-1)}(0,\lambda))_{j,\nu=1}^{n}=
\end{equation}%
\[
e^{inmt}+f_{1}(\lambda)e^{i(nm-1)t}+f_{2}(\lambda)e^{i(nm-2)t}+...+f_{nm-1}%
(\lambda)e^{it}+1
\]
which is a polynomial of $e^{it}$\ with entire coefficients $f_{1}%
(\lambda),f_{2}(\lambda),...$. Therefore the multiple eigenvalues of the
operators $L_{t}$ are the zeros of the resultant $R(\lambda)\equiv
R(\Delta,\Delta^{^{\prime}})$ of the polynomials $\Delta(\lambda,t)$ and
$\Delta^{^{\prime}}(\lambda,t)\equiv\frac{\partial}{\partial\lambda}%
\Delta(\lambda,t).$ Since $R(\lambda)$ is entire function and the large
eigenvalues of $L_{t}$ for $t\neq0,\pi$ are simple ( see Theorem 2 $(a)$),
\begin{equation}
\ker R=\{\lambda:R(\lambda)=0\}=\{a_{1},a_{2},...,\},\text{ }\lim
_{k\rightarrow\infty}a_{k}=\infty.
\end{equation}
For each $a_{k}$ there are $nm$ values $t_{k,1},t_{k,2},...,t_{k,nm}$ of $t$
satisfying $\Delta(a_{k},t)=0.$ Hence the set
\begin{equation}
A=\cup_{k=1}^{\infty}\{t:\Delta(a_{k},t)=0\}=\{t_{k,i}%
:i=1,2,...,nm;k=1,2,...,\}
\end{equation}
is countable and for $t\notin A$ all eigenvalues of $L_{t}$ are simple
eigenvalues. By Theorem 2$(a)$ the possible accumulation point of the set $A$
are $\pi k,$ where $k\in\mathbb{Z}.$

\begin{lemma}
The eigenvalues of $L_{t}$ can be numbered as $\lambda_{1}(t),$ $\lambda
_{2}(t),...,$ such that for each $p$ the function $\lambda_{p}(t)$ is
continuous in $Q$ and is analytic in $Q\backslash A(p),$ where $A(p)$ is a
subset of $A$ consisting of finite numbers $t_{1}^{p},t_{2}^{p},...,t_{s_{p}%
}^{p}.$ Moreover the followings hold:
\begin{equation}
\lim_{p\rightarrow\infty}\lambda_{p}(t)=\infty,\text{ }\lambda_{p(k,j)}%
(t)=\lambda_{k,j}(t),\text{ }\forall t\in Q_{\varepsilon}(n),
\end{equation}
where $|k|\geq N_{0},$ $p(k,j)=2|k|m+j$ if $k>0,$ $p(k,j)=(2|k|-1)m+j$ if
$k<0,$ the sets $Q_{\varepsilon}(n)$, $Q$ and number $N_{0}$ are defined in
(2) and in Theorem 1$(a).$
\end{lemma}

\begin{proof}
Let $t\in Q.$ It easily follows from the classical investigations [12, chapter
3, theorem 2] ( see (3), (4)) that there exist a large numbers $r$ and $c,$
independent of $t,$ such that the all eigenvalues of the operators $L_{t,z}$
for $z\in\lbrack0,1],$ where $L_{t,z}$ is defined by (45), lie in the set
\[
U(0,r)\cup(%
{\textstyle\bigcup\limits_{k:|k|\geq N_{0}}}
U(\left(  2\pi ki+ti\right)  ^{n},ck^{n-1-\frac{1}{2m}})),
\]
where $U(\mu,c)=\{\lambda\in\mathbb{C}:\mid\lambda-\mu\mid<c\}.$ Clearly there
exist a closed curve $\Gamma$ such that:

$(a)$ The curve $\Gamma$ lies in the resolvent set of the operators $L_{t,z}$
for all $z\in\lbrack0,1]$.

$(b)$ All eigenvalues of $L_{t,z}$ for all $z\in\lbrack0,1]$ that do not lie
in $U(\left(  2\pi ki+ti\right)  ^{n},ck^{n-1-\frac{1}{2m}})$ for $|k|\geq
N_{0}$ belong to the set enclosed by $\Gamma.$

Therefore taking into account that the family $L_{t,z}$ is holomorphic with
respect to $z,$ we obtain that the number of eigenvalues of operators
$L_{t,0}=L_{t}(C)$ and $L_{t,1}=L_{t}$ lying inside of $\Gamma$ are the same.
It means that apart from the eigenvalues $\lambda_{k,j}(t),$ where $|k|\geq
N_{0},$ $j=1,2,...,m,$ there exists $(2N_{0}-1)m$ eigenvalues of the operator
$L_{t}.$ We define $\lambda_{p}(t)$ for $p>(2N_{0}-1)m$ and $t\in
Q_{\varepsilon}(n)$ by (63). Let us first prove that these eigenvalues, that
is, the eigenvalues $\lambda_{k,j}(t)$ for $|k|\geq N_{0}$ are analytic
functions on $Q_{\varepsilon}(n).$ By Theorem 2$(a)$ if $t_{0}\in
Q_{\varepsilon}(n)$ and $|k|\geq N_{0}$ then $\lambda_{k,j}(t_{0})$ is a
simple root of (60), that is, $\Delta(\lambda,t_{0})=0,$ and $\Delta
^{^{\prime}}(\lambda,t_{0})\neq0$ for $\lambda=\lambda_{k,j}(t_{0}).$ By
implicit function theorem there exists a neighborhood $U(t_{0})$ of $t_{0}%
$\ and an analytic function $\lambda(t)$ on $U(t_{0})$ such that
$\Delta(\lambda(t),t)=0$ for $t\in U(t_{0})$ and $\lambda(t_{0})=\lambda
_{k,j}(t_{0}).$ By Theorem 2 $\lambda_{k,j}(t_{0})\in U(\mu_{k,j}(t_{0}%
),c_{1}|k|^{n-3}\ln|k|)).$ Since $\mu_{k,j}(t)$ and $\lambda(t)$ are
continuous functions the neighborhood $U(t_{0})$ of $t_{0}$ can be chosen so that

$\lambda(t)\in U(\mu_{k,j}(t),c_{1}|k|^{n-3}\ln|k|)$ for all $t\in U(t_{0}%
).$\ On the other hand, by Proposition 2, there exist unique eigenvalue of
$L_{t}$ lying in $U(\mu_{k,j}(t),c_{1}|k|^{n-3}\ln|k|)$ and this eigenvalue is
denoted by $\lambda_{k,j}(t).$ Therefore $\lambda(t)=\lambda_{k,j}(t)$ for all
$t\in U(t_{0})$, that is, $\lambda_{k,j}(t)$ is an analytic function in
$U(t_{0})$ for any $t_{0}\in Q_{\varepsilon}(n).$

Now let us continue analytically the function $\lambda_{p(k,j)}(t)$ into the
sets $U(0,\varepsilon),$ $U(\pi,\varepsilon)$ by using (60) and the implicit
function theorem. Consider (60) for
\[
t\in U(0,\varepsilon),\text{ }\lambda\in U_{0}=U((2\pi ki)^{n},2n(2\pi
k)^{n-1}\varepsilon).
\]
Since $U_{0}$ is a bounded region $(\ker R)\cap U_{0}$ is a finite set ( see
(61)). Therefore the subset $A(U_{0})$ of $A$ corresponding to $(\ker R)\cap
U_{0},$ that is, the values of $t$ corresponding to the multiple zeros of (60)
lying in $U_{0}$ \ is finite. It follows from (3) and (4) that for any $t\in
U(0,\varepsilon)\backslash A(U_{0})$ in the region $U_{0}$ the equation
$\Delta(\lambda,t)=0$ has$\ 2m$ different solutions $d_{1}(t),$ $d_{2}%
(t),...,d_{2m}(t)$ and
\[
\Delta^{^{\prime}}(\lambda,t)\neq0\text{ for }\lambda=d_{1}(t),d_{2}%
(t),...,d_{2m}(t).
\]
Using the implicit function theorem and taking into account (4) we see that
there exists a neighborhood $U(t,\delta)$ of $t$ such that:

$(i)$ There exist analytic functions $d_{1,t}(z),d_{2,t}(z),...,d_{2m,t}(z)$
in $U(t,\delta)$ coinciding with $d_{1}(t),d_{2}(t),...,d_{2m}(t)$ for $z=t$
respectively and satisfying%
\[
\Delta(d_{s,t}(z),z)=0,\text{ }d_{s,t}(z)\neq d_{j,t}(z),\forall z\in
U(t,\delta),\text{ }s=1,2,...,2m,\text{ }j\neq s.
\]

$(ii)$ $U(t,\delta)\cap A(U_{0})=\emptyset$ and $d_{s,t}(z)\in U_{0}$ for
$z\in U(t,\delta),$ $s=1,2,...,2m.$

Now take any point $t_{0}$ from $U(0,\varepsilon)\backslash A(U_{0})$. Let
$\gamma$ be line segment in $U(0,\varepsilon)\backslash A(U_{0})$ joining
$t_{0}$ and a point of the circle $S(0,\varepsilon)=\{t:|t|=\varepsilon\}.$
For any $t$ from $\gamma$ there exist $U(t,\delta)$ satisfying $(i)$ and
$(ii).$ Since $\gamma$ is a compact set the cover $\{U(t,\delta):t\in\gamma\}$
of $\gamma$ contains a finite cover $U(t_{0},\delta),U(t_{1},\delta
),...,U(t_{v},\delta),$ where $t_{v}\in S(0,\varepsilon).$ Now we are ready to
continue analytically the function $\lambda_{p(k,j)}(t)$ into the set
$U(0,\varepsilon).$ For any $z\in U(t_{v},\delta)\cap Q_{\varepsilon}(n)$ the
eigenvalue $\lambda_{p(k,j)}(z)$ coincides with one of the eigenvalues
$d_{1,t_{v}}(z),$ $d_{2,t_{v}}(z),...,$ $d_{2m,t_{v}}(z),$ since there exists
$2m$ eigenvalue of $L_{z}$ lying in $U_{0}.$ Denote by $B_{s}$ the subset of
the set $U(t_{v},\delta)\cap Q_{\varepsilon}(n)$ for which the function
$\lambda_{p(k,j)}(z)$ coincides with $d_{s,t_{v}}(z).$ Since $d_{s,t}(z)\neq
d_{i,t}(z)$ for $s\neq i$ the sets $B_{1},B_{2},...,B_{2m}$ are pairwise
disjoint and the union of these sets is $U(t_{v},\delta)\cap Q_{\varepsilon
}(n).$ Therefore there exists index $s$ for which the set $B_{s}$ contains
accumulation point and hence $\lambda_{p(k,j)}(z)=d_{s,t_{v}}(z)$ for all
$z\in U(t_{v},\delta)\cap Q_{\varepsilon}(n).$ Thus $d_{s,t_{v}}(z)$ is
analytic continuation of $\lambda_{p(k,j)}(z)$ to $U(t_{v},\delta).$ In the
same way we get the analytic continuation of $\lambda_{p(k,j)}(z)$ to
$U(t_{v-1},\delta),U(t_{v-2},\delta),...,U(t_{0},\delta).$ Since $t_{0}$ is
arbitrary point of $U(0,\varepsilon)\backslash A(U_{0})$ we obtain the
analytic continuation of $\lambda_{p(k,j)}(z)$ to $U(0,\varepsilon)\backslash
A(U_{0}).$ The analytic continuation of $\lambda_{p(k,j)}(z)$ to
$U(\pi,\varepsilon)\backslash A(U_{\pi})$ can be obtained in the same way,
where $A(U_{\pi})$ can be defined as $A(U_{0}).$ Thus the function
$\lambda_{p(k,j)}(t)$ is analytic in $Q\backslash A(p),$ where $A(p)$ consist
of finite numbers $t_{1}^{p},t_{2}^{p},...,t_{s_{p}}^{p}.$ Since
$\Delta(\lambda,t)$ is continuos with respect $(\lambda,t),$ the function
$\lambda_{p(k,j)}(t)$ can be extended continuously to the set $Q.$

Now let us define the eigenvalues $\lambda_{p}(t)$ for $p\leq(2N_{1}-1)m,$
$t\in Q$ which are apart from the eigenvalues defined by (63). These
eigenvalues lies in a bounded set $B$ and by (61) the set $B\cap\ker R$ and
the subset $A(B)$ of $A$ corresponding to $B$ are finite. Take a point $a$
from the set $Q\backslash A.$ Denote the eigenvalues of $L_{a}$ in increasing
( of absolute value) order $\mid\lambda_{1}(a)\mid\leq\mid\lambda_{2}%
(a)\mid\leq...\leq\mid\lambda_{(2N_{1}-1)m}(a)\mid.$ If $\mid\lambda
_{p}(a)\mid=\mid\lambda_{p+1}(a)\mid$ then by $\lambda_{p}(a)$ we denote the
eigenvalue that has a smaller argument, where argument is taken in $[0,2\pi)$.
Since $a\notin A$ the eigenvalues $\lambda_{1}(a),\lambda_{2}(a),...,\lambda
_{(2N_{1}-1)m}(a)$ are simple zeros of $\Delta(\lambda,a)=0.$ Therefore using
the implicit function theorem we obtain the analytic functions $\lambda
_{1}(t),\lambda_{2}(t),...,\lambda_{(2N_{1}-1)m}(t)$ on a neighborhood
$U(a,\delta)$ of $a$\ which are eigenvalues of $L_{t}$ for $t\in U(a,\delta).$
These functions can be analytically continued to $Q_{\varepsilon}(n)\backslash
A,$ being the eigenvalues of $L_{t},$ where, as we noted above, $A\cap
Q_{\varepsilon}(n)$ consist of a finite number of points. Taking into account
that $A(B)$ is finite, arguing as we have done in the proof of analytic
continuation and continuous extension of $\lambda_{p}(t)$ for $p>(2N_{1}-1)m$,
we obtain the analytic continuations of these functions to the set $Q$ except
finite points and continuous extension to $Q$
\end{proof}

By Gelfand's Lemma ( see [4]) every compactly supported vector function $f(x)$
can be represented in the form
\begin{equation}
f(x)=\frac{1}{2\pi}\int_{0}^{2\pi}f_{t}(x)dt,
\end{equation}
where $f_{t}(x)$ is defined by (53). This representation can be extended to
all function of \ $L_{2}^{m}(-\infty,\infty),$ and
\[
\int_{0}^{1}\langle f_{t}(x),X_{k,t}(x)\rangle dx=\int_{-\infty}^{\infty
}\langle f(x),X_{k,t}(x)\rangle dx,
\]
where $\{X_{k,t}:k=1,2,...\}$ is the biorthogonal system of $\{\Psi
_{k,t}:k=1,2,...\},$ $\Psi_{k,t}(x)$ is a normalized eigenfunction
corresponding to $\lambda_{k}(t),$ the eigenvalue $\lambda_{k}(t)$ is defined
in Lemma 6, $\Psi_{k,t}(x)$ and $X_{k,t}(x)$ are extended to $(-\infty
,\infty)$ by (58) and by $X_{k,t}(x+1)=e^{i\overline{t}}X_{k,t}(x)$.

Let $a\in(0,\frac{\pi}{2})\backslash A,$ $\varepsilon\in(0,\frac{a}{2})$ and
let $l(\varepsilon)$ be a smooth curve joining the points $-a$ and $2\pi-a$
and satisfying%
\begin{equation}
l(\varepsilon)\subset(Q_{\varepsilon}(n)\cap\Pi(a,\varepsilon))\backslash
A,\text{ }l(-\varepsilon)\cap A=\emptyset,\text{ }D(\varepsilon)\cup
\overline{D(-\varepsilon)}\subset Q
\end{equation}
where $\Pi(a,\varepsilon)=\{x+iy:x\in\lbrack-a,2\pi-a],y\in\lbrack
0,2\varepsilon)\},$ $l(-\varepsilon)=\{t:\bar{t}\in l(\varepsilon)\},$ the
sets $Q,$ $Q_{\varepsilon}(n)$ and $A$ are defined in (2) and (62),
$D(\varepsilon)$ and $D(-\varepsilon)$ are the domains enclosed by
$l(\varepsilon)\cup\lbrack-a,2\pi-a]$ and $l(-\varepsilon)\cup\lbrack
-a,2\pi-a]$ respectively, $\overline{D(-\varepsilon)}$ is closure of
$D(-\varepsilon).$ It is clear that, the domain $D(\varepsilon)\cup
\overline{D(-\varepsilon)}$ is enclosed by the closed curve $l(\varepsilon
)\cup l^{-}(-\varepsilon),$ where $l^{-}(-\varepsilon)$ is the opposite arc of
$l(-\varepsilon)$. Suppose $f\in S,$ that is, (11) holds. If $2\varepsilon
<\alpha$ then $f_{t}(x)$ is an analytic function of $t$ in a neighborhood of
$D(\varepsilon).$ Hence the Cauchy's theorem and (64) give
\begin{equation}
f(x)=\frac{1}{2\pi}\int_{l(\varepsilon)}f_{t}(x)dt.
\end{equation}
Since $l(\varepsilon)\in\mathbb{C}(n)$ ( see (65) and the definition of
$\mathbb{C}(n)$ in the introduction), it follows from Theorem 2($b$) and Lemma
6 that for each $t\in l(\varepsilon)$ we have a decomposition
\begin{equation}
f_{t}(x)=\sum_{k=1}^{\infty}a_{k}(t)\Psi_{k,t}(x),
\end{equation}
where $a_{k}(t)=(f_{t},X_{k,t}).$ Using (67) in (66) we get
\begin{equation}
f(x)=\frac{1}{2\pi}\int_{l(\varepsilon)}f_{t}(x)dt=\frac{1}{2\pi}%
\int_{l(\varepsilon)}\sum_{k=1}^{\infty}a_{k}(t)\Psi_{k,t}(x)dt.
\end{equation}

\begin{remark}
If $\lambda\in\sigma(L)$ \ then there exists points $t_{1},t_{2},...,t_{k}$ of
$[0,2\pi)$ such that $\lambda$ is an eigenvalue $\lambda(t_{j})$ \ of
$L_{t_{j}}$ of multiplicity $s_{j}$ for $j=1,2,...,k.$ Let $S(\lambda
,b)=\{z:\mid z-\lambda\mid=b\}$ be a circle containing only the eigenvalue
$\lambda(t_{j})$ \ of $L_{t_{j}}$ for $j=1,2,...,k.$ Using Lemma 6 we see that
there exists a neighborhood $U(t_{j},\delta)=\{t:\mid t-t_{j}\mid\leq\delta\}$
of $t_{j}$ such that:

$(a)$ The circle $S(\lambda,b)$ lies in the resolvent set of $L_{t}$ for all
$t\in U(t_{j},\delta)$ and $j=1,2,...,k.$

$(b)$ If $t\in(U(t_{j},\delta)\backslash\{t_{j}\}),$ then the operator $L_{t}$
has only $s_{j}$ eigenvalues lying in interior of $S(\lambda,b).$ These
eigenvalues are simple and let us denote they by $\Lambda_{j,1}(t),$
$\Lambda_{j,2}(t),...,\Lambda_{j,s_{j}}(t),$ where $j=1,2,...,k$.

Thus the spectrum of $L_{t}$ for $t\in U(t_{j},\delta),$ $j=1,2,...,k$
separated by $S(\lambda,b)$ into two parts in since of [7] ( see \S 6.4 of
chapter 3 of [7]). Since $\{L_{t}:$ $t\in U(t_{j},\delta)\}$ is a holomorphic
family of operators in since [7] (see \S 1 of chapter 7 of [7]), the theory of
holomorphic family of finite dimensional operators can be applied to the part
of $L_{t}$ for $t\in U(t_{j},\delta)$ corresponding to the inside of
$S(\lambda,b).$ Therefore ( see \ \S 1 of the chapter 2 of [7] ) the
eigenvalue $\Lambda_{j,1}(t),$ $\Lambda_{j,2}(t),...,$ $\Lambda_{j,s_{j}}(t)$
and corresponding eigenprojections $P(\Lambda_{j,1}(t)),$ $P(\Lambda
_{j,2}(t)),...,P(\Lambda_{j,s_{j}}(t))$ are branches of an analytic function.
These eigenprojections is represented by a Laurent series in $t^{\frac{1}{\nu
}},$ where $\nu\leq s_{j},$ with finite principal parts. One can easily see
that if $\lambda_{p}(t)$ is a simple eigenvalue of $L_{t}$ then
\begin{equation}
P(\lambda_{p}(t))f=(f,X_{p,t})\Psi_{p,t},\text{ }\parallel P(\lambda
_{p}(t))\parallel=\frac{1}{\parallel X_{p,t}\parallel}=\mid\frac{1}{\alpha
_{p}(t)}\mid
\end{equation}
and $P(\lambda_{p}(t))$ is analytic function in some neighborhood of $t,$
where $\alpha_{p}(t)=(\Psi_{p,t},\Psi_{p,t}^{\ast}).$ This and Lemma 6 show
that for each $p$ the function $a_{p}(t)\Psi_{p,t}$ is analytic on
$D(\varepsilon)\cup\overline{D(-\varepsilon)}$ except finite points.
\end{remark}

\begin{theorem}
$(a)$ If $f(x)$ is absolutely continuous, compactly supported function and

$f^{^{\prime}}\in L_{2}^{m}(-\infty,\infty)$ then
\begin{equation}
f(x)=\frac{1}{2\pi}\sum_{k=1}^{\infty}\int_{l(\varepsilon)}a_{k}(t)\Psi
_{k,t}(x)dt
\end{equation}
and
\begin{equation}
f(x)=\frac{1}{2\pi}\sum_{k=1}^{\infty}\int_{[0,2\pi)^{+}}a_{k}(t)\Psi
_{k,t}(x)dt,
\end{equation}
where
\begin{equation}
\int_{\lbrack0,2\pi)^{+}}a_{k}(t)\Psi_{k,t}(x)dt=\lim_{\varepsilon
\rightarrow0}\int_{l(\varepsilon)}a_{k}(t)\Psi_{k,t}(x)dt.
\end{equation}
and the series (70), (71) converge uniformly in any bounded subset of
$(-\infty,\infty).$

$(b)$ Every function $f(x)\in S,$ where $S$ is defined in (11), has
decompositions (70) and (71), where the series converges in the norm of
$L_{2}^{m}(a,b)$ for every $a,b\in\mathbb{R}.$
\end{theorem}

\begin{proof}
The proof of (70) in case $(a)$ follows from (68), Theorem 3, and Lemma 6. In
Appendix A by writing the proof of \ the Theorem 2 of [19] in the vector form
we get the proof of (70) in the case $(b).$ In Appendix B the formula (71) is
obtained from (70) by writing the proof of \ the Theorem 3 of [19] in the
vector form
\end{proof}

\begin{definition}
Let $\lambda$ be a point of the spectrum $\sigma(L)$ of $L$ and $t_{1}%
,t_{2},...,t_{k}$ be the points of $[0,2\pi)$ such that $\lambda$ is a
eigenvalue of $L_{t_{j}}$ of multiplicity $s_{j}$ for $j=1,2,...,k.$ The point
$\lambda$ is called a spectral singularity of $L$ if
\begin{equation}
\sup\parallel P(\Lambda_{j,i}(t))\parallel=\infty,
\end{equation}
where supremum is taken over all $t\in(U(t_{j},\delta)\backslash\{t_{j}\}),$
$j=1,2,...,k;$ $i=1,2,...,s_{j},$ the set $U(t_{j},\delta)$ and the
eigenvalues $\Lambda_{j,1}(t),$\ $\Lambda_{j,2}(t),...,$ $\Lambda_{j,s_{j}%
}(t)$ are defined in Remark 1. In other words $\lambda$ is called a spectral
singularity of $L$ if there exists indices $j,i$ such that the point $t_{j}$
is a pole of $P(\Lambda_{j,i}(t)).$ Briefly speaking a point $\lambda\in
\sigma(L)$ is called a spectral singularity of $L$ if the projections of
$L_{t}$ corresponding to the simple eigenvalues lying in the small
neighborhood of $\lambda$ are not uniformly bounded. We denote the set of
spectral singularities by $S(L).$
\end{definition}

\begin{remark}
Note that if $\gamma=\{\lambda_{p}(t):t\in(\alpha,\beta)\}$ is a curve lying
in $\sigma(L)$ and containing no multiple eigenvalues of $L_{t},$ where
$t\in\lbrack0,2\pi),$ then arguing as in papers [18,9] one can prove that for
the projection $P(\gamma)$ of $L$ corresponding to $\gamma$ the following
hold
\begin{equation}
P(\gamma)f=\int_{(\alpha,\beta)}(f,X_{p,t})\Psi_{p,t},dt,\text{ }\parallel
P(\gamma)\parallel=\sup_{t\in(\alpha,\beta)}\parallel P(\lambda_{p}%
(t))\parallel,
\end{equation}
that is, the definition 1 is equivalent to the definition of the spectral
singularities given in [18,9], where the spectral singularities is defined as
a points in the neighborhoods of which the projections $P(\gamma)$ are not
uniformly bounded. The proof of (74) is long technical. In order to avoid
eclipsing the essence by technical detail and taking into account that in the
spectral expansion of $L$ the eigenfunctions and eigenprojections of $L_{t}$
for $t\in\lbrack0,2\pi)$ are used ( see (71)), and using that there are the
closed relationship between projections (see (74)) of $L$ and $L_{t}$ for
$t\in\lbrack0,2\pi),$ in this paper, in the definition of the spectral
singularities, without loss of naturalness, instead of the boundlessness of
projections $P(\gamma)$ of $L$ we use the boundlessness of projections
$P(\lambda_{p}(t))$, of $L_{t},$ that is, we use the definition 1. \ In any
case the spectral singularity is a point of $\sigma(L)$ that requires the
regularization in order to get the spectral expansion.
\end{remark}

\begin{theorem}
$(a)$ All spectral singularity of $L$ are contained in the set of the multiply
eigenvalues of $L_{t}$ for $t\in\lbrack0,2\pi)$, that is, $S(L)=\{\Lambda
_{1},\Lambda_{2},...\}\subset\ker R\cap\sigma(L),$ where $S(L)$ and $\ker R$
are defined in the Definition 1 and in (61).

$(b)$ Let $\lambda=\lambda_{p}(t_{0})\in\sigma(L)\backslash S(L),$ where
$t_{0}\in(a,2\pi-a).$ If $\gamma_{1},$ $\gamma_{2},...,$ are sequence of
smooth curves lying in a neighborhood $U=\{t\in\mathbb{C}$: $\mid t-t_{0}%
\mid\leq\delta_{0}\}$ of $t_{0}$ and approximating the interval $[t_{0}%
-\delta_{0},t_{0}+\delta_{0}]$ then
\begin{equation}
\lim_{k\rightarrow\infty}%
{\textstyle\int\limits_{\gamma_{k}}}
a_{p}(t)\Psi_{p,t}(x)dt=%
{\textstyle\int\limits_{t_{0}-\delta_{0}}^{t_{0}-\delta_{0}}}
a_{p}(t)\Psi_{p,t}(x)dt,
\end{equation}
where $U$ is a neighborhood of $t_{0}$ such that if $t\in$ $U$ then
$\lambda_{p}(t)$ is not a spectral singularity.

$(c)$ If the operator $L$ has not spectral singularities then we have the
following spectral expansion in term of the parameter $t:$
\begin{equation}
f(x)=\frac{1}{2\pi}\sum_{k=1}^{\infty}\int_{0}^{2\pi}a_{k}(t)\Psi_{k,t}(x)dt.
\end{equation}
If $f(x)$ is absolutely continuous, compactly supported function and
$f^{^{\prime}}\in L_{2}^{m}(-\infty,\infty)$ then the series in (76) converges
uniformly in any bounded subset of $(-\infty,\infty).$ If $f(x)\in S$ then the
series converges in the norm of $L_{2}^{m}(a,b)$ for every $a,b\in\mathbb{R}.$
\end{theorem}

\begin{proof}
$(a)$ If $\lambda_{p}(t_{0})$ is a simple eigenvalue of $L_{t_{0}}$ then due
to the Remark 1 ( see (69) and the end of Remark 1) the projection
$P(\lambda_{p}(t))$ and $\mid\alpha_{p}(t)\mid$ continuously depend on $t$ in
some neighborhood of $t_{0}.$ On the other hand $\alpha_{p}(t_{0})\neq0,$
since the system of the root functions of $L_{t_{0}}$ is complete. Therefore
it follows from the Definition 1 that $\lambda$ is not a spectral
singularities of $L.$

$(b)$ It follows from (61) and Theorem 5$(a)$ that there exists a neighborhood
$U$ of $t_{0}$ such that if $t\in$ $U$ then $\lambda_{p}(t)$ is not spectral
spectral singularities of $L.$ If $\lambda_{p}(t_{0})\in\sigma(L)\backslash
S(L)$ then by Definition 1 $t_{0}$ is not a pole of $P(\lambda_{p}(t)),$ that
is, by Remark 1 the Laurent series in $t^{\frac{1}{\nu}},$ where $\nu\leq s,$
of $P(\lambda_{p}(t))$ at $t_{0}$ has not principal part. Therefore (69)
implies that $\frac{1}{\mid\alpha_{p}(t)\mid}$ and hence $\frac{1}{\mid
\alpha_{p}(t)\mid}(f_{t},\Psi_{p,t}^{\ast})\Psi_{p,t}$ is a bounded continuous
functions in some neighborhood of $t_{0}$ , which implies the proof of $(b).$

$(c)$ It follows from Theorem 5$(b)$ that if the operator $L$ has not spectral
singularities then
\begin{equation}
\int_{\lbrack0,2\pi)^{+}}a_{k}(t)\Psi_{k,t}(x)dt=\int_{0}^{2\pi}a_{k}%
(t)\Psi_{k,t}(x)dt,
\end{equation}
where the left-hand side is defined by (72). Thus (76) follows from (77), (71)
\end{proof}

Now we change the variables to $\lambda$ by using the characteristic equation
$\Delta(\lambda,t)=0$ and the implicit-function theorem. By (60)
$\Delta(\lambda,t)$ and $\frac{\partial\Delta(\lambda,t)}{\partial t}$ are
polynomials of $e^{it}$ and their resultant is entire function. It is clear
that this resultant is not zero function. Let $b_{1},b_{2},...,$ be zeros of
the resultant, i.e., are the common zeros of the polynomials $\Delta
(\lambda,t)$ and $\frac{\partial\Delta(\lambda,t)}{\partial t}.$ Then
$\lim_{k\rightarrow\infty}b_{k}=\infty$ and the equation $\Delta(\lambda,t)=0$
defines a function $t(\lambda)$ such that
\begin{equation}
\Delta(\lambda,t(\lambda))=0,\text{ }\frac{dt}{d\lambda}=-\frac{\partial
\Delta/\partial\lambda}{\partial\Delta/\partial t},\text{ }\frac
{\partial\Delta(\lambda,t)}{\partial t}/_{t=t(\lambda)}\neq0\text{ }%
\end{equation}
for all $\lambda\in\mathbb{C}\backslash\{b_{1},b_{2},...,\}.$ Consider the
functions
\begin{equation}
F_{p,t}(x)=\sum_{k=1,2...,n}Y_{k}(x,\lambda_{p}(t))A_{k}(t,\lambda
_{p}(t))=(\sum_{k=1,2...,n}Y_{k}(x,\lambda)A_{k}(t(\lambda),\lambda
))_{\lambda=\lambda_{p}(t)}%
\end{equation}
where $Y_{1}(x,\lambda),Y_{2}(x,\lambda),\ldots,Y_{n}(x,\lambda)$ are linearly
independent solutions of (59),

$A_{k}=(A_{k,1},A_{k,2},...,A_{k,m}),$ $A_{k,i}=A_{k,i}(t,\lambda)$ is the
cofactor of the entry in $mn$ row and $(k-1)m+i$ column of the determinant
(60). One can readily see that%
\begin{equation}
A_{k,i}(t,\lambda)=g_{s}(\lambda)e^{ist}+g_{s-1}(\lambda)e^{i(s-1)t}%
+...+g_{1}(\lambda)e^{it}+g_{0}(\lambda),
\end{equation}
where $g_{0}(\lambda),g_{1}(\lambda),...,$ are entire functions. By (78)
$A_{k,i}(t(\lambda),\lambda)$ is analytic function in $\mathbb{C}%
\backslash\{b_{1},b_{2},...,\}.$ Since the operator $L_{t}$ for $t\neq0,\pi$
has a simple eigenvalue there exists a nonzero cofactor of the determinant
(60). Without loss of generality it can be assumed that $A_{k,1}%
(t(\lambda),\lambda)$ is nonzero function. Then $A_{k,1}(t(\lambda),\lambda)$
has a finite number zeros in each compact subset of $\mathbb{C}\backslash
\{b_{1},b_{2},...,\}.$\ Therefore there exists a countable set $E_{1}$ such
that
\begin{equation}
\{b_{1},b_{2},...,\}\subset E_{1},\text{ }A_{k,1}(t(\lambda),\lambda
)\neq0,\text{ }\forall\lambda\notin E_{1}.
\end{equation}
Let $A_{1}$ be the set of all $t$ satisfying $\Delta(\lambda,t)=0$ for some
$\lambda\in E_{1}.$ Clearly $A_{1}$ is a countable set. Now using Lemma 6,
(79), (81) and taking into account that the functions $Y_{1}(x,\lambda
),Y_{2}(x,\lambda),\ldots,Y_{n}(x,\lambda)$ are linearly independent, we
obtain
\begin{equation}
\Psi_{p,t}(x)=\frac{F_{p,t}(x)}{\parallel F_{p,t}\parallel},\text{ }\parallel
F_{p,t}\parallel\neq0,\text{ }\forall t\in(D(\varepsilon)\cup\overline
{D(-\varepsilon)})\backslash(A\cup A_{1}),
\end{equation}
where $\Psi_{p,t}(x)$ is a normalized eigenfunction corresponding to
$\lambda_{p}(t).$ Since the set $A\cup A_{1}$ is countable there exist the
curves $l(\varepsilon_{1}),l(\varepsilon_{2}),...,$ such that
\begin{equation}
\lim_{s\rightarrow\infty}l(\varepsilon_{s})=[-a,2\pi-a],\text{ }%
l(\varepsilon_{s})\in(D(\varepsilon)\cup\overline{D(-\varepsilon)}%
)\backslash(A\cup A_{1}),\text{ }\forall s.
\end{equation}
Now let us do the change of variables in (70). Using (78), (79), (82) we get
\[
a_{p}(t(\lambda))\Psi_{p,t(\lambda)}(x)=\frac{h(\lambda)}{\alpha(\lambda
)}F(x,\lambda),
\]
where $F(x,\lambda)=\sum_{j=1,2...,n}Y_{j}(x,\lambda)A_{j}(\lambda),$
$A_{j}(\lambda)=A_{j}(t(\lambda),\lambda),$ ( see (79), (80) for the
definition of $A_{j}(t,\lambda)),$ $(F(x,\lambda))_{\lambda=\lambda_{p}%
(t)}=F_{p,t}(x),$ $h(\lambda)=(f(\cdot),\Phi(\cdot,\lambda))$, $\Phi
(x,\lambda_{p}(t))$ is eigenfunction of $L_{t}^{\ast}$ corresponding to
$\overline{\lambda_{p}(t)\text{ }}$ and $\alpha(\lambda)\equiv(F(\cdot
,\lambda),\Phi(\cdot,\lambda)).$ Using this notations and (78), we obtain
\begin{equation}%
{\displaystyle\int\limits_{l(\varepsilon_{s})}}
a_{p}(t)\Psi_{p,t}(x)dt=%
{\displaystyle\int\limits_{\Gamma_{p}(\varepsilon_{s})}}
\frac{-h(\lambda)\varphi(\lambda)}{\alpha(\lambda)\phi(\lambda)}(\sum
_{j=1}^{n}Y_{j}(x,\lambda)A_{j}(\lambda))d\lambda,
\end{equation}
where $\Gamma_{p}(\varepsilon_{s})=\{\lambda=\lambda_{p}(t):t\in
l(\varepsilon_{s})\},$ $\varphi=\partial\Delta/\partial\lambda,$
$\phi=\partial\Delta/\partial t$ . Note that it follows from (78) and (83)
that $\phi(\lambda)\neq0$ for $\lambda\in\Gamma_{p}(\varepsilon_{s}).$ If $t$
$\in l(\varepsilon_{s})$ then by the definition of $A$ and by (83)
$\lambda_{p}(t)$ is a simple eigenvalue. Hence $\alpha_{p}(t)\neq0,$ since the
root functions of $L_{t}$ is complete in $L_{2}^{m}(0,1).$ Therefore
$\alpha(\lambda)\neq0$ for $\lambda\in\Gamma_{p}(\varepsilon_{s}).$

To do the regularization about the spectral singularities $\Lambda_{1}%
,\Lambda_{2},...,$ we take into account that there are numbers $i_{l}$ and
$\delta$ such that for $\mid\lambda-\Lambda_{l}\mid<\delta$ the equality
\[
\mid\frac{(\lambda-\Lambda_{l})^{i_{l}}h(\lambda)\varphi(\lambda)A_{j}%
(\lambda)}{\alpha(\lambda)\phi(\lambda)}\mid<c_{18}%
\]
for $j=1,2,...,n$ holds and the neighborhoods $U_{\delta}(\Lambda
_{l})=\{\lambda:\mid\lambda-\Lambda_{l}\mid<\delta\}$ do not intersect.
Introduce the mapping $B$ as follows:
\[
Bf(x,\lambda)=f(x,\lambda)-\sum_{l}\sum_{\nu=0}^{i_{l}-1}B_{l,\nu}%
(\lambda)\frac{\partial^{\nu}(f(x,\Lambda_{l}))}{\partial\lambda^{\nu}},
\]
where $B_{l,\nu}(\lambda)=\frac{(\lambda-\Lambda_{l})^{\nu}}{\nu!}$ for
$\lambda\in U_{\delta,}(\Lambda_{l})$ and $B_{l,\nu}(\lambda)=0$ for
$\lambda\notin U_{\delta}(\Lambda_{l}).$ We set
\[
\Gamma_{k}=\{\lambda=\lambda_{k}(t):t\in\lbrack0,2\pi)\},\text{ }%
S_{k}=\{l:\Lambda_{l}\in\Gamma_{k}\cap S(L)\}.
\]
Now using this notations and formulas (71), (72), (84), we get
\begin{equation}
f(x)=\frac{1}{2\pi}\sum_{k=1}^{\infty}(\int_{\Gamma_{k}}\frac{-h(\lambda
)\varphi(\lambda)}{\alpha(\lambda))\phi(\lambda)}(\sum_{j=1}^{n}%
B(Y_{j}(x,\lambda))A_{j}(\lambda))d\lambda+\sum_{l\in S_{k}}M_{k,l}(x)),
\end{equation}
where
\[
M_{k,l}(x)=\lim_{s\rightarrow\infty}\frac{1}{2\pi}\int_{\Gamma_{k}%
(\varepsilon_{s})}\frac{-h(\lambda)\varphi(\lambda)}{\alpha(\lambda
))\phi(\lambda)}(\sum_{j=1}^{n}(\sum_{\nu=0}^{i_{l}-1}B_{l,\nu}(\lambda
)\frac{\partial^{\nu}(Y_{j}(x,\Lambda_{l}))}{\partial\lambda^{\nu}}%
A_{j}(\lambda))d\lambda.
\]
Thus Theorem 4 implies the following spectral expansion of $L:$

\begin{theorem}
Every function $f(x)\in S$ has decomposition (85), where the series converges
in the norm of $L_{2}^{m}(a,b)$ for every $a,b\in\mathbb{R}.$ If $f(x)$ is
absolutely continuous, compactly supported function and $f^{^{\prime}}\in
L_{2}^{m}(-\infty,\infty)$ then the series in (85) converges uniformly in any
bounded subset of $(-\infty,\infty).$
\end{theorem}

\begin{remark}
If $n=2\mu+1,$ then by Theorem 2 all large eigenvalue of $L_{t}$ for $t\in Q$
are simple, the set $A\cap Q$, is finite, the number of spectral singulariries
is finite ( if exist), (77) holds for $k\gg1,$ and if $\varepsilon$ is small
number, then $D_{\varepsilon}$ and $D_{\varepsilon}^{-}$ do not contain the
point of $A$. Therefore the spectral expansion (85) has a simpler form.
Moreover repeating the proof of \ Corollary 1$(a)$ of [21], we obtain that
every function $f(x)$ satisfying (10) has decomposition (85).
\end{remark}

\section{Appendices}

APPENDIX \ A. THE PROOF OF (70)

Here we justify the term by term integration of the series in (68). Let
$H_{N,t}$ be the linear span of $\Psi_{1,t}(x),$ $\Psi_{2,t}(x),...,\Psi
_{N,t}(x)$ and $f_{N,t}$ be the projection of $f_{t}(x)$ onto $H_{N,t}.$ Since
$\{\Psi_{k,t}(x)\}$ and $\{X_{k,t}(x)\}$ are biorthogonal system we have
\begin{equation}
f_{N,t}(x)=\sum_{k=1,2,...,N}a_{k}^{N}(t)\Psi_{k,t}(x), \tag{A1}%
\end{equation}
where $a_{k}^{N}(t)=(f_{N,t},X_{k,t}).$ Using the notations $g_{N,t}%
=f_{t}-f_{N,t},$ $b_{k}^{N}(t)=(g_{N,t},X_{k,t})$ and (A1), we obtain
$a_{k}^{N}(t)=a_{k}(t)-b_{k}^{N}(t)$ and
\[
f_{t}=\sum_{k=1,2,...,N}(a_{k}(t)-b_{k}^{N}(t))\Psi_{k,t}+g_{N,t}.
\]
This with (66) give
\begin{align}
f(x)  &  =\frac{1}{2\pi}(\sum_{k=1,2,...,N}\int_{l_{\varepsilon}}(a_{k}%
(t)\Psi_{k,t}(x)dt\tag{A2}\\
&  +\int_{l_{\varepsilon}}(g_{N,t}(x)-\sum_{k=1,2,..,N}b_{k}^{N}(t))\Psi
_{k,t}(x))dt).\nonumber
\end{align}
To obtain (70) we must to prove that the last integral in (A2) tends to zero
as $N\rightarrow\infty.$ For this we prove the following

\begin{lemma}
On $l_{\varepsilon}$ the functions
\begin{equation}
\parallel g_{N,t}\parallel,\text{ }\parallel\sum_{k=1,2,...,N}b_{k}%
^{N}(t))\Psi_{k,t}\parallel\tag{A3}%
\end{equation}
tend to zero as $N\rightarrow\infty$ uniformly with respect to $t.$
\end{lemma}

\begin{proof}
First we prove that $\parallel g_{N,t}\parallel$ tends to zero uniformly. Let
$P_{N,t}$ and $P_{\infty,t}$ be projections of $L_{2}^{m}[0,1]$ onto $H_{N,t}$
and $H_{\infty,t}$ respectively, where $H_{\infty,t}=\cup_{n=1}^{\infty
}H_{N,t}.$ If follows from (67) that $f_{t}\in H_{\infty,t}.$ On the other
hand one can readily see that
\begin{equation}
H_{N,t}\subset H_{N+1,t}\subset H_{\infty,t},\text{ }P_{N,t}\subset
P_{\infty,t},\text{ }P_{N,t}\rightarrow P_{\infty,t}. \tag{A4}%
\end{equation}
Therefore $P_{N,t}f_{t}\rightarrow f_{t},$ that is $\parallel g_{N,t}%
\parallel\rightarrow0.$ Since $\parallel g_{N,t}\parallel$ is a distance from
$f_{t}$ to $H_{N,t},$ for each sequence $\{t_{1},t_{1},...\}\subset
l(\varepsilon)$ converging to $t_{0}$ we have
\begin{align*}
&  \parallel g_{N,t_{s}}\parallel\leq\parallel f_{t_{s}}-\sum_{k=1,2,...,N}%
a_{k}^{N}(t_{0})\Psi_{k,t_{s}}(x)\parallel\leq\parallel g_{N,t_{0}}%
\parallel+\\
&  \parallel f_{t_{s}}-f_{t_{0}}\parallel+\parallel\sum_{k=1,2,...,N}a_{k}%
^{N}(t_{0})(\Psi_{k,t_{0}}-\Psi_{k,t_{s}})\parallel\leq\parallel g_{N,t_{0}%
}\parallel+\alpha_{s},
\end{align*}
where $\alpha_{s}\rightarrow0$ as $s\rightarrow\infty$ by continuity of
$f_{t}$ and $\Psi_{k,t}$ on $l(\varepsilon).$ Similarly ( interchanging
$t_{0}$ and $t_{s}),$ we get $\parallel g_{N,t_{0}}\parallel\leq\parallel
g_{N,t_{s}}\parallel+\beta_{s},$ where $\beta_{s}\rightarrow0$ as
$s\rightarrow\infty.$ Hence $\parallel g_{N,t}\parallel$ is a continuos
function on the compact $l(\varepsilon).$ On the other hand the first
inclusion \ of (A4) implies that $\parallel g_{N,t}\parallel\geq\parallel
g_{N+1,t}\parallel.$ \ Now it follows from the proved three properties of
$\parallel g_{N,t}\parallel$ that $\parallel g_{N,t}\parallel$ tend to zero as
$N\rightarrow\infty$ uniformly on the compact $l(\varepsilon).$

Now to prove that the second function in (A3) tends to zero uniformly we
consider the family of operators $\Gamma_{p,t}$ for $t\in l(\varepsilon),$
$p=1,2,...,$ by formula
\begin{equation}
\Gamma_{p,t}(f)=\sum_{k=1,2,...,p}(f,X_{k,t})\Psi_{k,t}(x). \tag{A5}%
\end{equation}
First let us prove that the set
\begin{equation}
\Gamma(f)=\{\Gamma_{p,t}(f):t\in l(\varepsilon),p=1,2,...,\} \tag{A6}%
\end{equation}
is a bounded subset of $L_{2}^{m}[0,1].$ Since in the Hilbert space every
weakly bounded subset is a strongly bounded subset, it is enough to show that
for each $g\in L_{2}^{m}[0,1]$ there exists a constant $M$ such that
\begin{equation}
\mid(g,\varphi)\mid<M,\forall\varphi\in\Gamma(f). \tag{A7}%
\end{equation}
Decomposing $g$ by the basis $\{X_{k,t}:k=1,2,...,\}$, using definition of
$\varphi$ ( see (A7), (A6), (A5)), and then the uniform asymptotic formulas
(8), (9) we obtain%
\[
\mid(g,\varphi)\mid\leq\sum_{k=1,2,...,p}\mid(\varphi,X_{k,t})(g,\Psi
_{k,t})\mid\leq\sum_{k=1,2,...,p}\mid(\varphi,X_{k,t})\mid^{2}%
\]%
\[
+\sum_{k=1,2,...,p}\mid(g,\Psi_{k,t})\mid^{2}\leq\parallel\varphi\parallel
^{2}+\parallel g\parallel^{2}+c_{19}.
\]
which implies (A7). Thus $\Gamma(f)$ is a bounded set. On the other hand one
can readily see that $\Gamma_{p,t}$ for $t\in l(\varepsilon),$ $p=1,2,...,$ is
a linear continuous operator. Therefore by Banach- Steinhaus theorem the
family of operators $\Gamma_{p,t}$ is equicontinuous. Now using the equality
\[
\Gamma_{N,t}g_{N,t}=\sum_{k=1,2,...,N}b_{k,j}^{N}(t))\Psi_{k,j,t}%
\]
and taking into account that the first function in (A3) tend to zero
uniformly, we obtain that the second function in (A3) also tends to zero uniformly
\end{proof}

Using Lemma 7 and Schwarz inequality we get%
\[
\parallel\int_{l_{\varepsilon}}(g_{N,t}(x)-\sum_{k=1,2,...,N}b_{k}^{N}%
(t))\Psi_{k,t}(x))dt\parallel\leq
\]%
\[
C_{\varepsilon}\int_{a}^{b}\int_{l(\varepsilon)}\mid g_{N,t}(x)-\sum
_{k=1,2,...,N}b_{k,}^{N}(t))\Psi_{k,t}(x))\mid\mid dt\mid dx=
\]%
\[
C_{\varepsilon}\int_{l(\varepsilon)}\parallel(g_{N,t}(x)-\sum_{k=1,2,...,N}%
b_{k}^{N}(t))\Psi_{k,t}(x))\parallel\mid dt\mid\rightarrow0\text{ as
}N\rightarrow\infty,
\]
where $C_{\varepsilon}$ is the length of $l(\varepsilon),$ the norm used here
is the norm of $L_{2}^{m}(a,b),$ $a$ and $b$ are the real numbers. This with
(A2) justify the term by term integration of the series in (68).

APPENDIX B. THE PROOF OF (71)

Here we use the notation introduced in (65) and prove (71). Since for fixed
$k$ the function $a_{k}(t)\Psi_{k,t}(x)$ is analytic on $D(\varepsilon)$
except finite number points $t_{1}^{k},t_{2}^{k},...,t_{p_{k}}^{k}$ ( see the
end of the Remark 1) we have
\begin{equation}
\int_{l(\varepsilon)}a_{k}(t)\Psi_{k,t}dt=\int_{[0,2\pi]^{+}}a_{k}%
(t)\Psi_{k,t}dt+\sum_{s:t_{s}^{k}\in D(\varepsilon)}\text{Res}_{t=t_{s}^{k}%
}a_{k}(t)\Psi_{k,t}, \tag{B1}%
\end{equation}
Similarly
\begin{equation}
\int_{l(-\varepsilon)}a_{k}(t)\Psi_{k,t}dt=\int_{[0,2\pi]^{+}}a_{k}%
(t)\Psi_{k,t}dt+\sum_{s:t_{s}^{k}\in\overline{D(-\varepsilon)}}\text{Res}%
_{t=t_{s}^{k}}a_{k}(t)\Psi_{k,t}. \tag{B2}%
\end{equation}
Since $l(\varepsilon)\cup l^{-}(-\varepsilon)$ is a closed curve enclosing
$D(-\varepsilon)\cup\overline{D(-\varepsilon)},$ we have
\begin{equation}
\int_{l(\varepsilon)\cup l^{-}(-\varepsilon)}a_{k,j}(t)\Psi_{k,t}%
(x)dt=\sum_{s:t_{s}^{k}\in D(-\varepsilon)\cup\overline{D(-\varepsilon)}%
}\text{Res}_{t=t_{s}^{k}}a_{k}(t)\Psi_{k,t}. \tag{B3}%
\end{equation}
Now applying (70) to the curves $l(\varepsilon),$ $l(-\varepsilon),$
$l(\varepsilon)\cup l^{-}(-\varepsilon)$, using (B1), (B2), (B3) and taking
into account that $l(\varepsilon)\cup l^{-}(-\varepsilon)$ is a closed curve,
we obtain
\begin{equation}
f(x)=\frac{1}{2\pi}\sum_{k=1,2,...}(\int_{[0,2\pi]^{+}}a_{k}(t)\Psi
_{k,t}(x)dt+\sum_{s:t_{s}^{k}\in D(\varepsilon)}\text{Res}_{t=t_{s}^{k}}%
a_{k}(t)\Psi_{k,t}), \tag{B4}%
\end{equation}%
\begin{equation}
f(x)=\frac{1}{2\pi}\sum_{k=1,2,...}(\int_{[0,2\pi]^{+}}a_{k}(t)\Psi
_{k,t}(x)dt+\sum_{s:t_{s}^{k}\in\overline{D(-\varepsilon)}}\text{Res}%
_{t=t_{s}^{k}}a_{k}(t)\Psi_{k,t}). \tag{B5}%
\end{equation}%
\begin{equation}
0=\frac{1}{2\pi}\int_{l(\varepsilon)\cup l^{-}(-\varepsilon)}f_{t}%
(x)dt=\frac{1}{2\pi}\sum_{k=1,2,...}(\sum_{s:t_{s}^{k}\in(D(-\varepsilon
)\cup\overline{D(-\varepsilon)})}\text{Res}_{t=t_{s}^{k}}a_{k}(t)\Psi_{k,t}).
\tag{B6}%
\end{equation}
Adding (B4) and (B5) and then using (B6) we get the proof of (71).

\end{document}